\documentclass[reqno,final,11pt]{amsart}
\usepackage{natbib}  
\usepackage{fancyhdr} 
\usepackage{color} 
\usepackage{hyperref} 
\usepackage{graphicx} 
\usepackage{geometry}
\usepackage[utf8]{inputenc}
\usepackage[T1]{fontenc}
\usepackage{verbatim}

\usepackage{tikz}
\usetikzlibrary{calc,shapes,backgrounds,arrows}
\usepackage{amssymb}
\usepackage{enumitem}
\usepackage{stmaryrd}

\definecolor{aleacolor}{rgb}{0.16,0.59,0.78}

\hypersetup{
breaklinks,
colorlinks=true,
linkcolor=aleacolor,
urlcolor=aleacolor,
citecolor=aleacolor}

\renewcommand{\cite}{\citet}

\theoremstyle{plain}
\newtheorem{theorem}{Theorem}[section]                                          
\newtheorem{proposition}[theorem]{Proposition}                          
\newtheorem{lemma}[theorem]{Lemma}
\newtheorem{corollary}[theorem]{Corollary}

\newtheorem{open}[theorem]{Open question}
\theoremstyle{definition}
\newtheorem{definition}[theorem]{Definition}
\theoremstyle{remark}
\newtheorem{remark}[theorem]{Remark}

\makeatletter \@addtoreset{equation}{section} \makeatother
\renewcommand\theequation{\thesection.\arabic{equation}}
\renewcommand\thefigure{\thesection.\arabic{figure}}
\renewcommand\thetable{\thesection.\arabic{table}}

\newcommand{\aleaIndex}[1]{\href{http://alea.impa.br/english/index_v18.htm}{\bf 18}}
\newcommand{\aleaDOI}[1]{\href{https://doi.org/10.30757/ALEA.v18-}{DOI: 10.30757/ALEA.v18-}}
\eheader{ALEA, Lat. Am. J. Probab. Math. Stat.}{\aleaIndex{18}}{2021}{1}{4}{\aleaDOI}
\elogo{\parbox[c]{180pt}{\includegraphics[width=170pt]{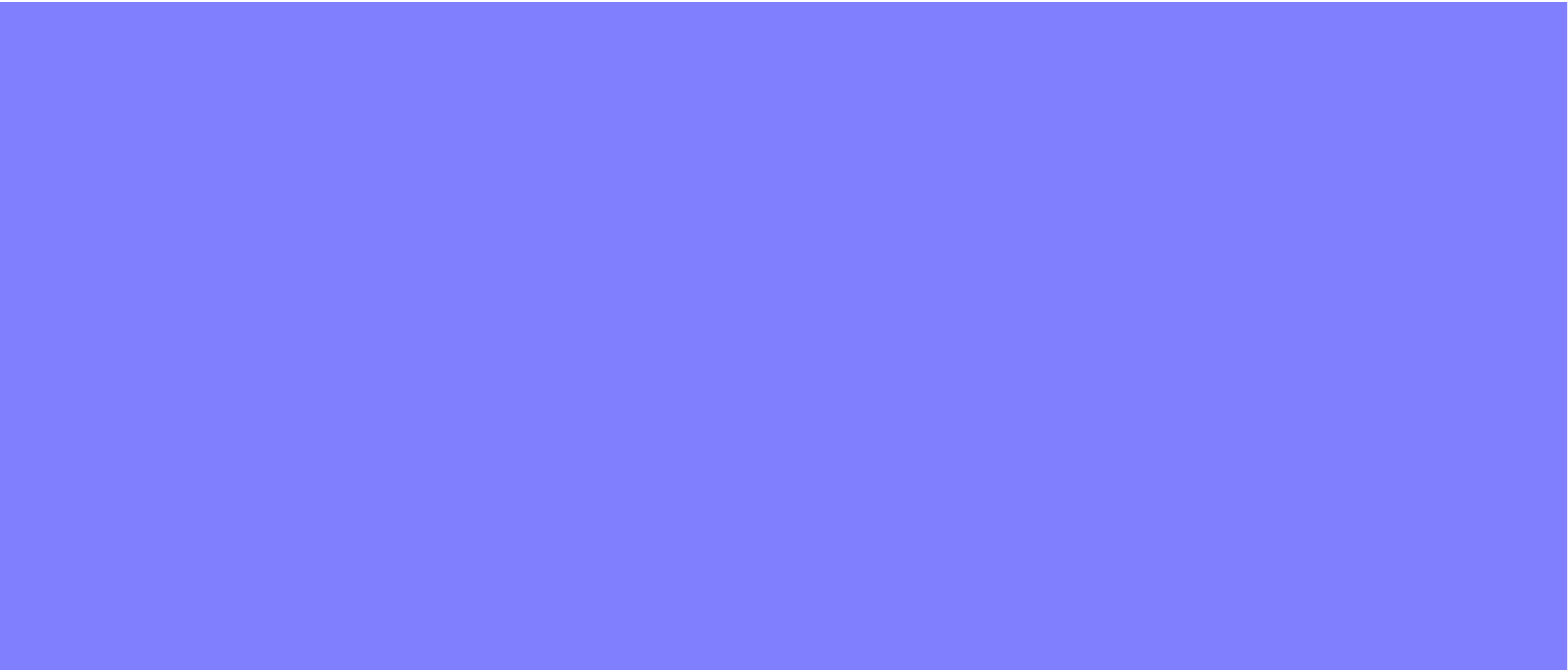}}} 

\begin{document}

\title[Fractal properties of the frontier in Poissonian coloring]{Fractal properties of the frontier in Poissonian coloring}

\author{Anne-Laure Basdevant}
\address{Laboratoire de Probabilités Statistique et Modélisation, Sorbonne Université, 4 place Jussieu, 75005 Paris, France}
\email{anne.laure.basdevant@normalesup.org} 
\urladdr{\href{http://https://www.normalesup.org/~basdevant/}{https://www.normalesup.org/$\sim$basdevant/}} 

\author{Guillaume Blanc}
\address{Laboratoire de Mathématiques d'Orsay, Université Paris-Saclay, rue Michel Magat, 91400 Orsay, France}
\email{guillaume.blanc1@universite-paris-saclay.fr} 
\urladdr{\href{https://sites.google.com/view/guillaume-blanc-math}{https://sites.google.com/view/guillaume-blanc-math}} 

\author{Nicolas Curien}
\address{Laboratoire de Mathématiques d'Orsay, Université Paris-Saclay, rue Michel Magat, 91400 Orsay, France}
\email{nicolas.curien@gmail.com} 
\urladdr{\href{https://www.imo.universite-paris-saclay.fr/~nicolas.curien/}{https://www.imo.universite-paris-saclay.fr/$\sim$nicolas.curien/}} 

\author{Arvind Singh}
\address{Laboratoire de Mathématiques d'Orsay, CNRS et Université Paris-Saclay, rue Michel Magat, 91400 Orsay, France}
\email{arvind.singh@universite-paris-saclay.fr} 
\urladdr{\href{https://www.imo.universite-paris-saclay.fr/~arvind.singh/}{https://www.imo.universite-paris-saclay.fr/$\sim$arvind.singh/}} 

\thanks{The first and fourth authors were supported by ANR 19-CE40-0025 ProGraM. 
The second and third authors were supported by ERC 740943 GeoBrown and ANR RanTanPlan.}

\begin{abstract}
We study a model of random partitioning  by nearest-neighbor coloring from Poisson rain, introduced independently by \cite{aldous2018random} and \cite{preater2009species}. 
Given two initial points in $[0,1]^{d}$ respectively colored in red and blue, we let independent uniformly random points fall in $[0,1]^d$, and upon arrival, each point takes the color of the nearest point fallen so far. 
We prove that the colored regions converge in the Hausdorff sense towards two random closed subsets whose intersection, the \emph{frontier}, has Hausdorff dimension strictly between $d-1$ and $d$, thus answering a conjecture raised by \cite{aldous2018random}. 
However, several topological properties of the frontier remain elusive.
\end{abstract}

\maketitle
 
\section{Introduction and main results}

\begin{figure}[ht]
\begin{center}
\includegraphics[height=3cm]{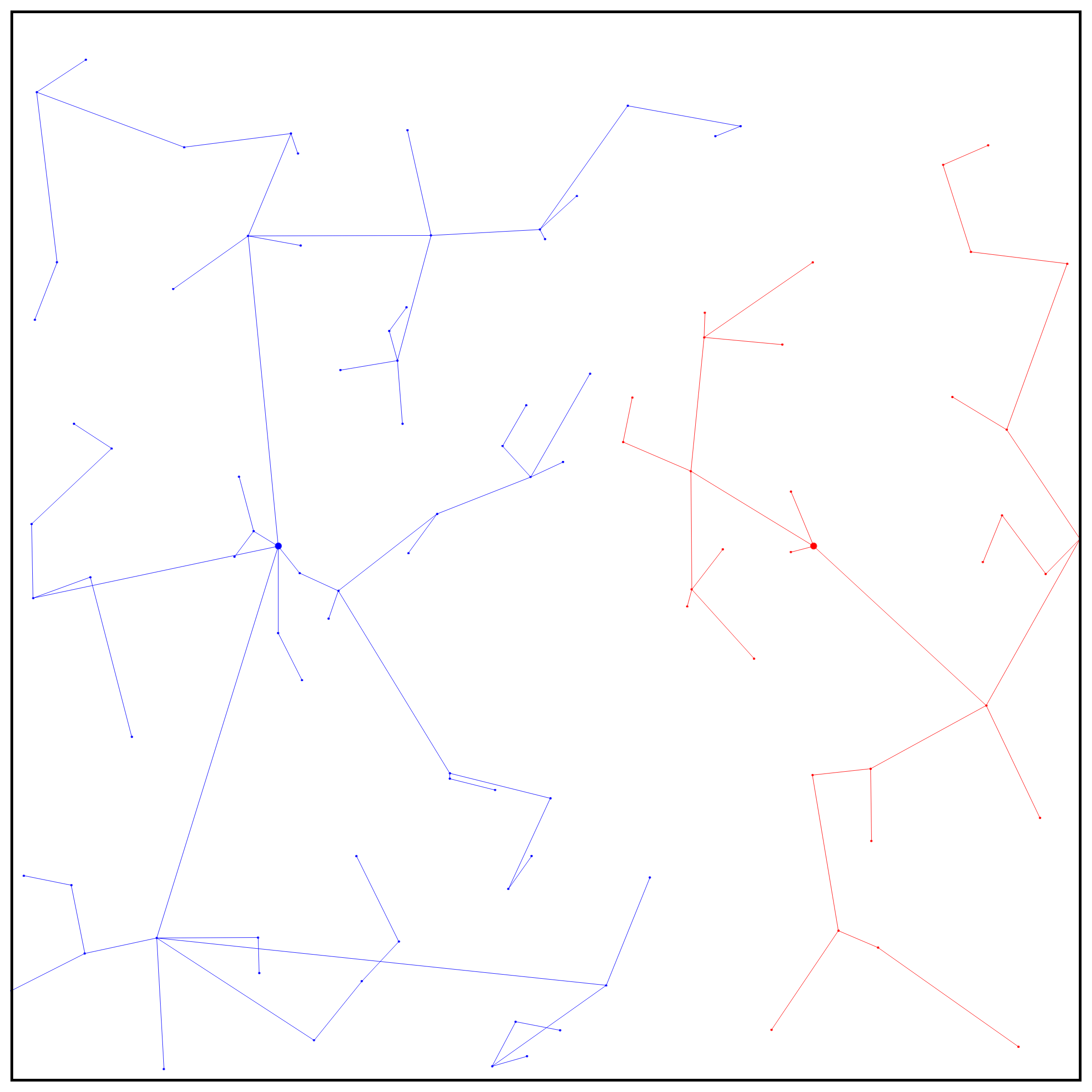}
\includegraphics[height=3cm]{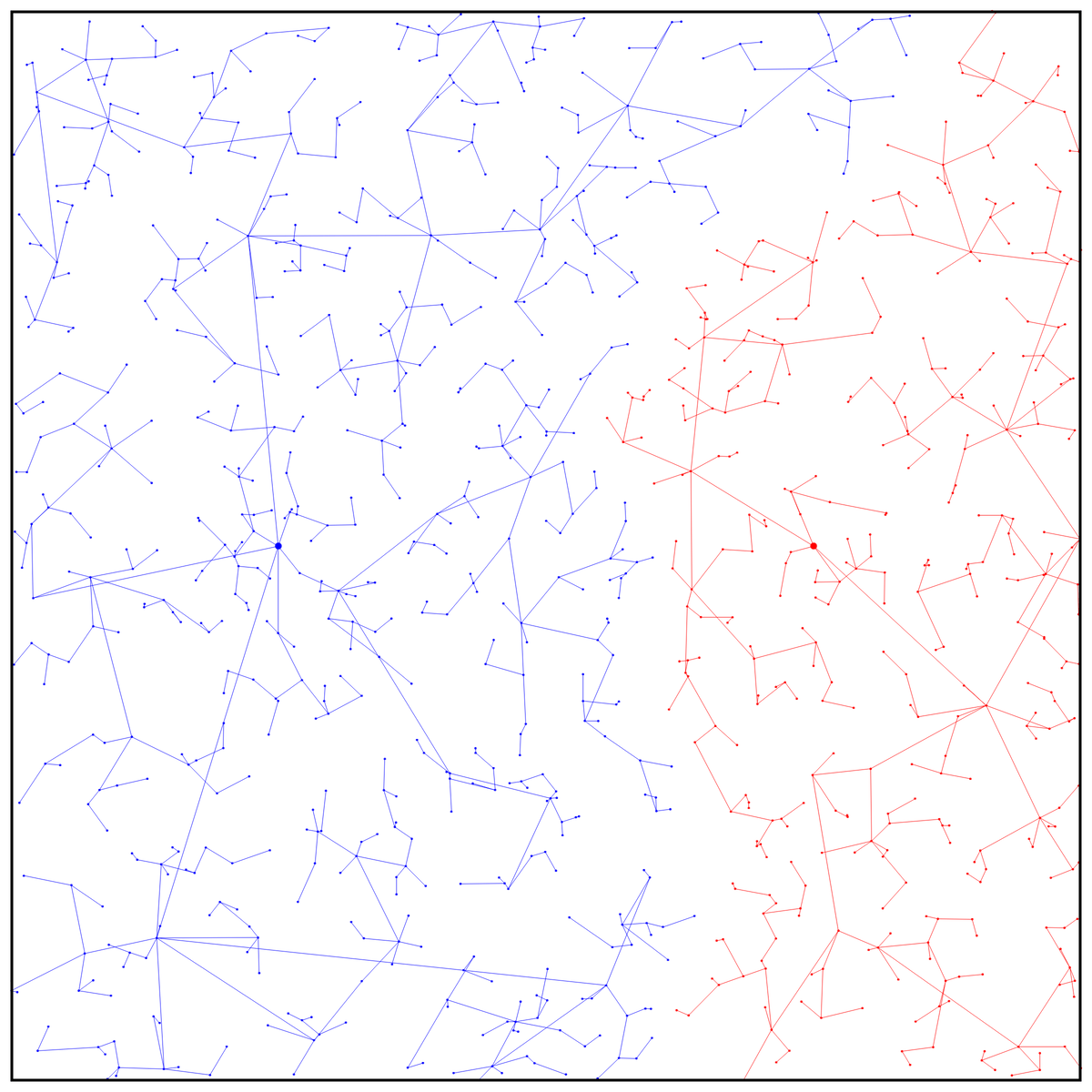}
\includegraphics[height=3cm]{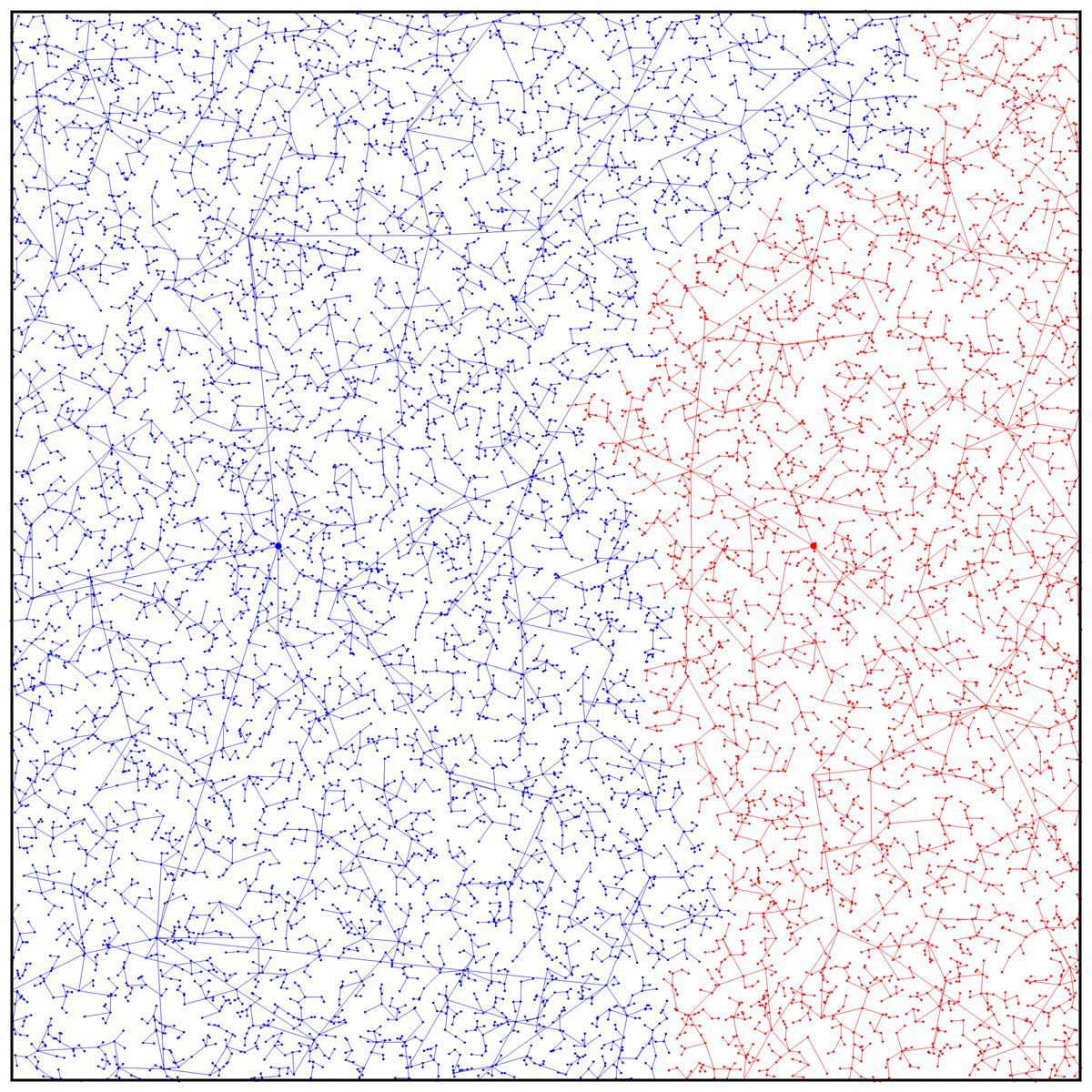}
\includegraphics[height=3cm]{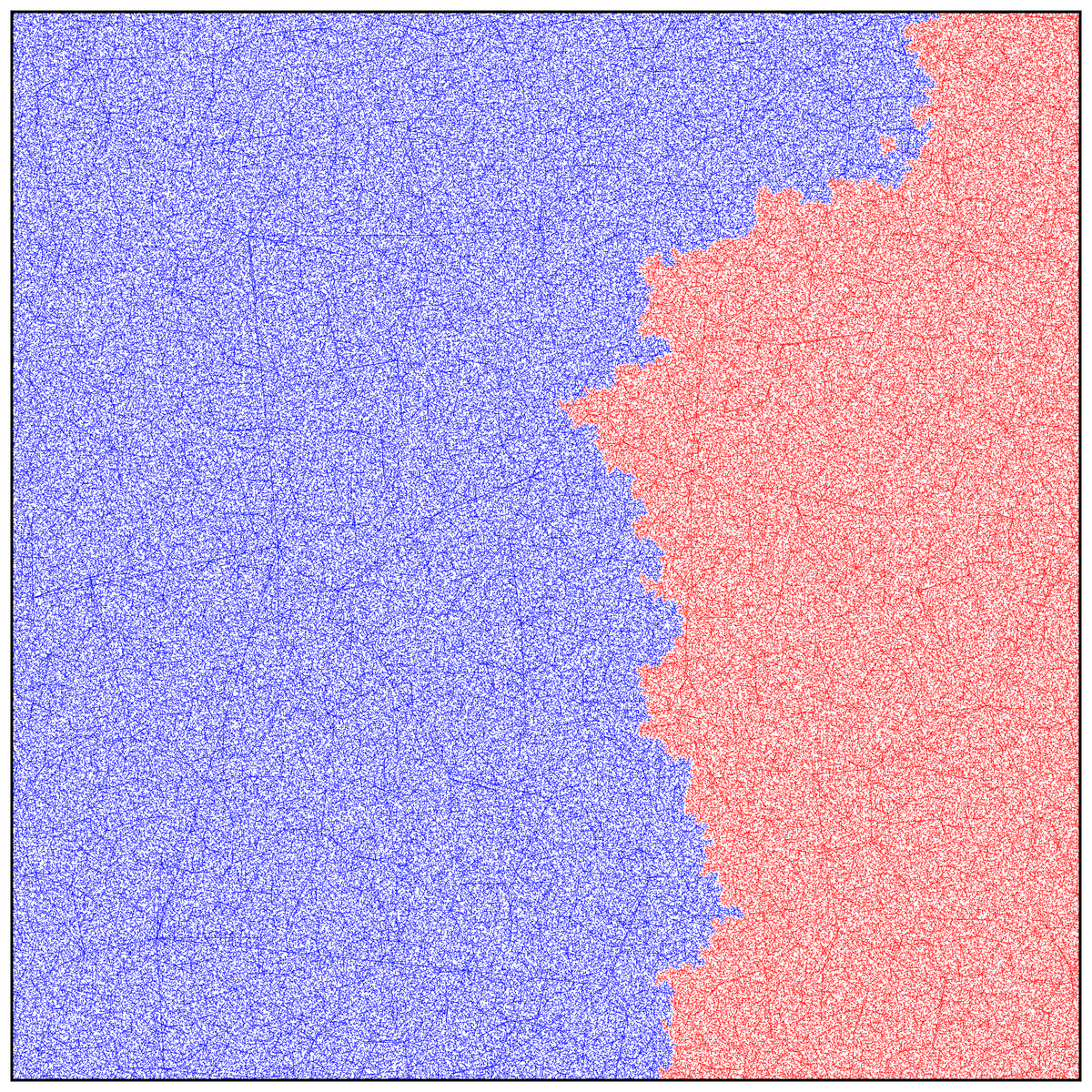}
\includegraphics[height=3cm]{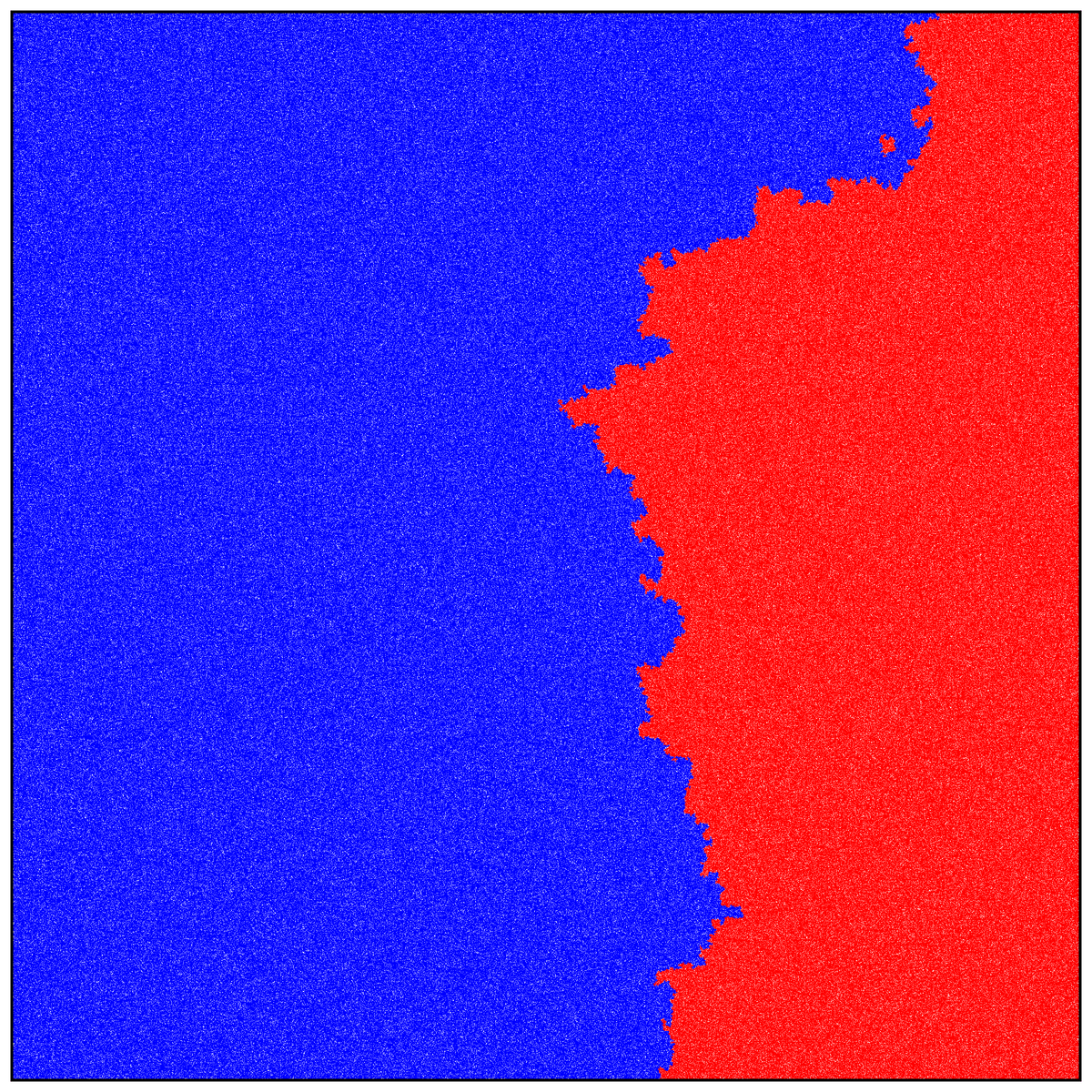}
\caption{Simulation of the Poisson coloring of space where a new incoming point takes the color of the nearest neighbor in the process so far, from left to right with $10^{2}, 10^{3},10^{4},10^{6}$ and $10^{7}$ points. \label{fig:simpoint}}
\end{center}
\end{figure}

We consider a model of Poissonian coloring which is based on a dynamical construction in the $d$-dimensional hypercube $[0,1]^d$. 
Initially, two points $R_0 \ne B_0$  are planted in $[0,1]^d$: think of $R_0$ as an initial red seed, and of $B_0$ as an initial blue seed. 
All the randomness in the construction comes from a sequence $(X_n)_{n\in\mathbb{N}^*}$ of independent random variables, uniformly distributed in $[0,1]^d$.
Picturing $X_1,X_2,\ldots$ as points falling consecutively in $[0,1]^d$, we let each point take the color of the closest point already present (nearest neighbor for the usual Euclidean metric $d$). 
Formally, define the initial red and blue sets as $\mathcal{R}_0=\{R_0\}$ and $\mathcal{B}_0=\{B_0\}$, respectively.
Then, by induction, for each $n\in\mathbb{N}$ such that the red and blue sets $\mathcal{R}_n$ and $\mathcal{B}_n$ have been constructed, proceed as follows: almost surely, we have $d(X_{n+1},\mathcal{R}_n)\neq d(X_{n+1},\mathcal{B}_n)$, and
\begin{itemize}
\item if $d(X_{n+1},\mathcal{R}_n)<d(X_{n+1},\mathcal{B}_n)$, then set $\mathcal{R}_{n+1}=\mathcal{R}_n\cup\{X_{n+1}\}$ and $\mathcal{B}_{n+1}=\mathcal{B}_n$;
\item otherwise, if $d(X_{n+1},\mathcal{R}_n)>d(X_{n+1},\mathcal{B}_n)$, then set $\mathcal{R}_{n+1}=\mathcal{R}_n$ and $\mathcal{B}_{n+1}=\mathcal{B}_n\cup\{X_{n+1}\}$.
\end{itemize}
Letting $n\to\infty$, the red and blue sets $\mathcal{R}_n$ and $\mathcal{B}_n$ respectively converge, for the Hausdorff distance between closed subsets of $[0,1]^d$, to
\[\mathcal{R}_\infty=\overline{\bigcup_{n\geq0}\mathcal{R}_n}\quad\text{and}\quad\mathcal{B}_\infty=\overline{\bigcup_{n\geq0}\mathcal{B}_n}.\]
We refer to \cite{burago^2ivanov} for background on the Hausdorff distance (in particular, the above convergence is implied by Exercise 7.3.5 therein).
The object we are interested in is the \emph{frontier} $\mathcal{F}_\infty=\mathcal{R}_\infty\cap\mathcal{B}_\infty$, which is also easily shown to be the limit as $n\to\infty$ of the discrete frontier 
\[\mathcal{F}_n=\left\{x\in[0,1]^d:d(x,\mathcal{R}_n)=d(x,\mathcal{B}_n)\right\},\]
still for the Hausdorff distance (c.f.~Proposition \ref{prop:Fnontrivial} below). 
See Figure \ref{fig:simpoint} for a simulation of the coloring process.

\medskip

This very natural model can be found in \cite{aldous2018random}, where it is attributed to \cite[Section 7.6.8]{penrose2009random}, although it may have been considered by other authors before.
Recently, \cite{LichevMitsche} studied the combinatorial properties of genealogical trees induced by the coloring procedure. 
Here, we focus instead on the geometric and topological properties of the model. 
After completion of this work, we learned from Aldous that \cite{preater2009species} had considered the same model, and in particular answered \cite[Conjecture 3]{aldous2018random}, showing that the frontier $ \mathcal{F}_\infty$ has zero Lebesgue measure (see \cite[Theorem 2]{preater2009species}). 
Our main result is the following, and settles a conjecture of \cite[Section 5.3.3]{aldous2018random}.

\begin{theorem}[The frontier is fractal]\label{thm:main} 
Almost surely, the Hausdorff dimension of the frontier $\mathcal{F}_\infty$ satisfies
$$d-1<\dim_H\mathcal{F}_{\infty}<d.$$
\end{theorem}

The proof is divided into two main steps, which we summarize below.

\begin{itemize}
\item \textsc{Upper bound}\textbf{.} We first show that for every $x\in[0,1]^d$ and $r>0$ such that the ball $\overline{B}(x,r)$ does not contain the seeds $R_0$ and $B_0$, there is a positive probability that the smaller ball $\overline{B}(x,r/6)$ is monochromatic at the end of the coloring (see Lemma \ref{lem:isolated} below). 
Together  with a multi-scale argument, this shows that the Hausdorff dimension of the frontier $\mathcal{F}_\infty$ is strictly less than $d$ (Proposition \ref{propkeylemma'}). 
A result with a similar flavor, also using a first-passage percolation argument, was obtained by Preater (see \cite[Theorem 1]{preater2009species}), who then showed that $\mathcal{F}_\infty$ has zero Lebesgue measure.
\item \textsc{Lower bound}\textbf{.} The lower bound on the Hausdorff dimension of the frontier is based on ideas and techniques developed in \cite[Sections 5 and 6]{aizenman1999holder}, where the authors introduce general conditions which allow to lower bound the Hausdorff dimension of random curves (see \cite[Theorem 1.3]{aizenman1999holder}). 
Their result applies in particular to scaling limits of interfaces from critical statistical physics models such as percolation; random curves which have a positive probability, at each scale, of oscillating. Unfortunately, it is not clear that our frontier $ \mathcal{F}_\infty$ even contains curves, see Open question \ref{open:curve} below. We find a workaround by adapting the ideas of Aizenman and Burchard to get a Hausdorff dimension lower bound result for connected random closed subsets of $[0,1]^d$. 
The exact statement is given in Theorem \ref{thmAB'}. 
We hope that this extension can prove to be of independent interest.
\end{itemize}

\medskip

\textbf{A natural variant.}
There are natural variants of this coloring model, such as the following ``segment'' model (as opposed to the original ``point'' model): still thinking of $R_0$ and $B_0$ as initial red and blue seeds, and of $X_1,X_2,\ldots$ as points falling consecutively in $[0,1]^d$, let as before $\mathcal{R}_0=\{R_0\}$ and $\mathcal{B}_0=\{B_0\}$ be the initial red and blue sets, respectively. Then, by induction, for each $n\in\mathbb{N}$ such that the red and blue sets $\mathcal{R}_n$ and $\mathcal{B}_n$ have been constructed, proceed as follows: almost surely, we have $d(X_{n+1},\mathcal{R}_n)\neq d(X_{n+1},\mathcal{B}_n)$, and
\begin{itemize}
\item if $d(X_{n+1},\mathcal{R}_n)<d(X_{n+1},\mathcal{B}_n)$, then set $\mathcal{R}_{n+1}=\mathcal{R}_n\cup[Y_n,X_{n+1}]$ and $\mathcal{B}_{n+1}=\mathcal{B}_n$, where $Y_n$ denotes the point on $\mathcal{R}_n$ which is closest to $X_{n+1}$;
\item otherwise, if $d(X_{n+1},\mathcal{R}_n)>d(X_{n+1},\mathcal{B}_n)$, then set $\mathcal{R}_{n+1}=\mathcal{R}_n$ and $\mathcal{B}_{n+1}=\mathcal{B}_n\cup[Y_n,X_{n+1}]$, where $Y_n$ denotes the point on $\mathcal{B}_n$ which is closest to $X_{n+1}$.
\end{itemize}
Note that, by construction, the red and blue sets $\mathcal{R}_n$ and $\mathcal{B}_n$ are connected finite unions of line segments, so that $Y_n$ is always well defined (such a point is almost surely unique because $X_{n+1}$ is uniform and independent of $X_1,\ldots,X_n$).
Upon minor technical modifications in the proofs, Theorem \ref{thm:main} holds for this coloring process as well.

\begin{figure}[ht]
\begin{center}
\includegraphics[height=3cm]{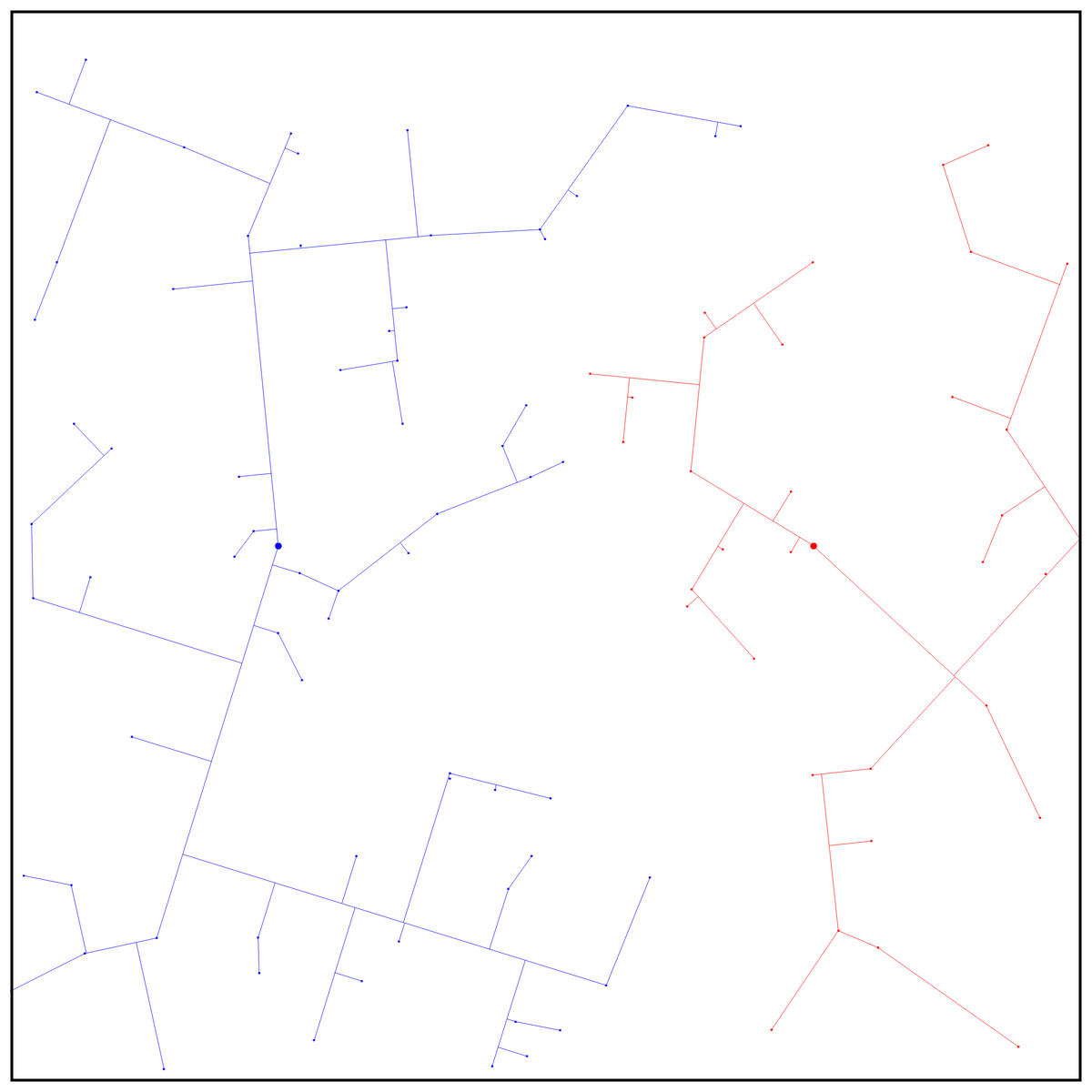}
\includegraphics[height=3cm]{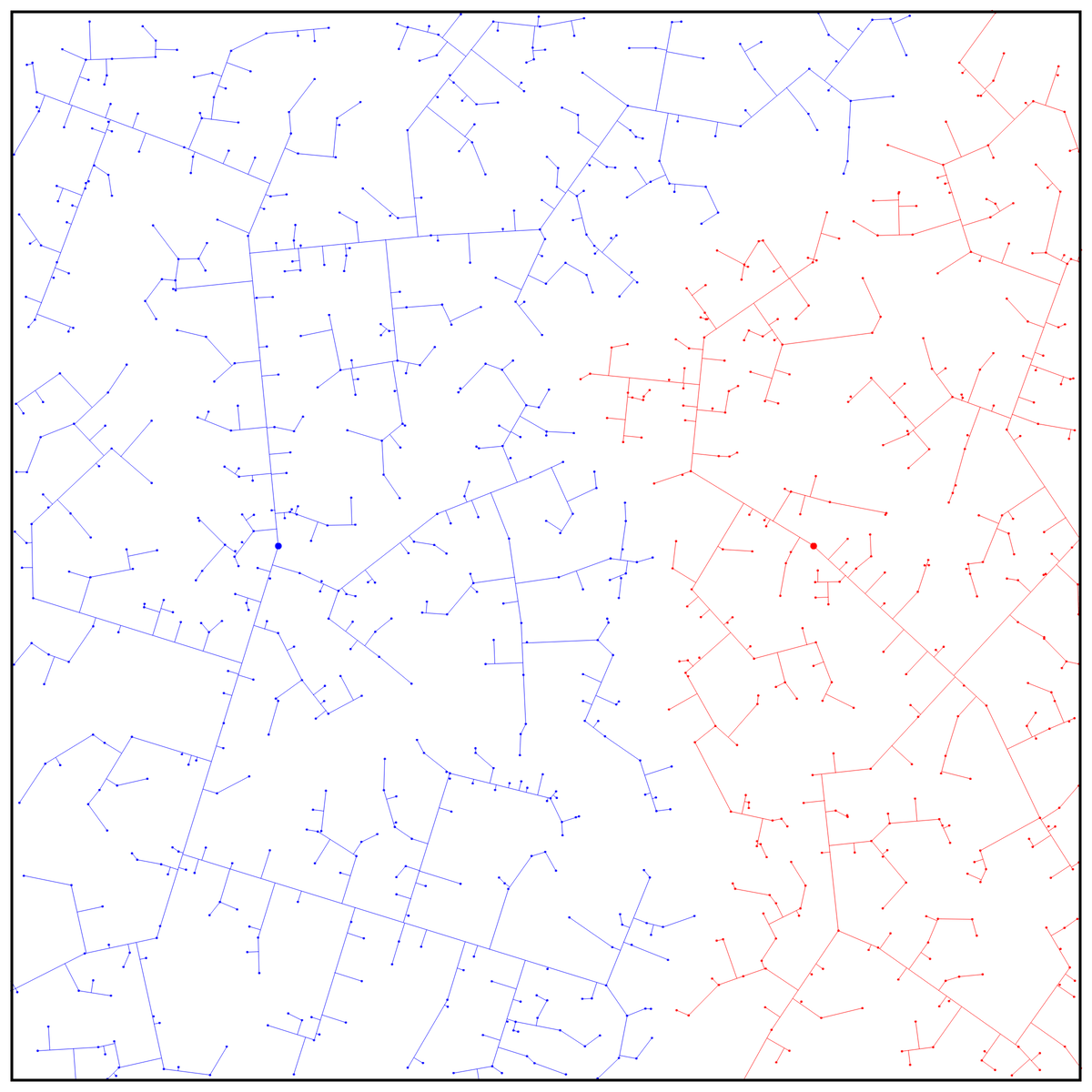}
\includegraphics[height=3cm]{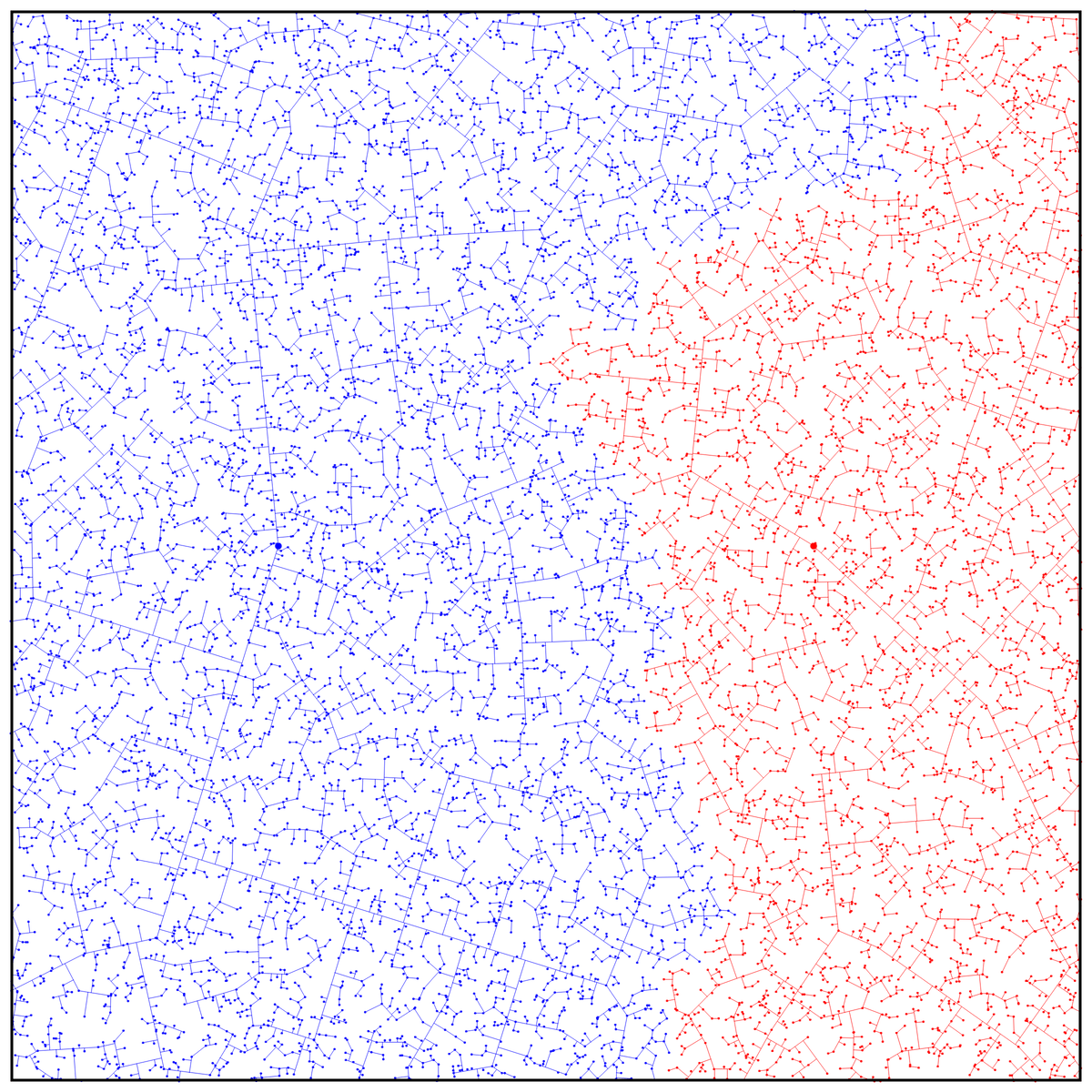}
\includegraphics[height=3cm]{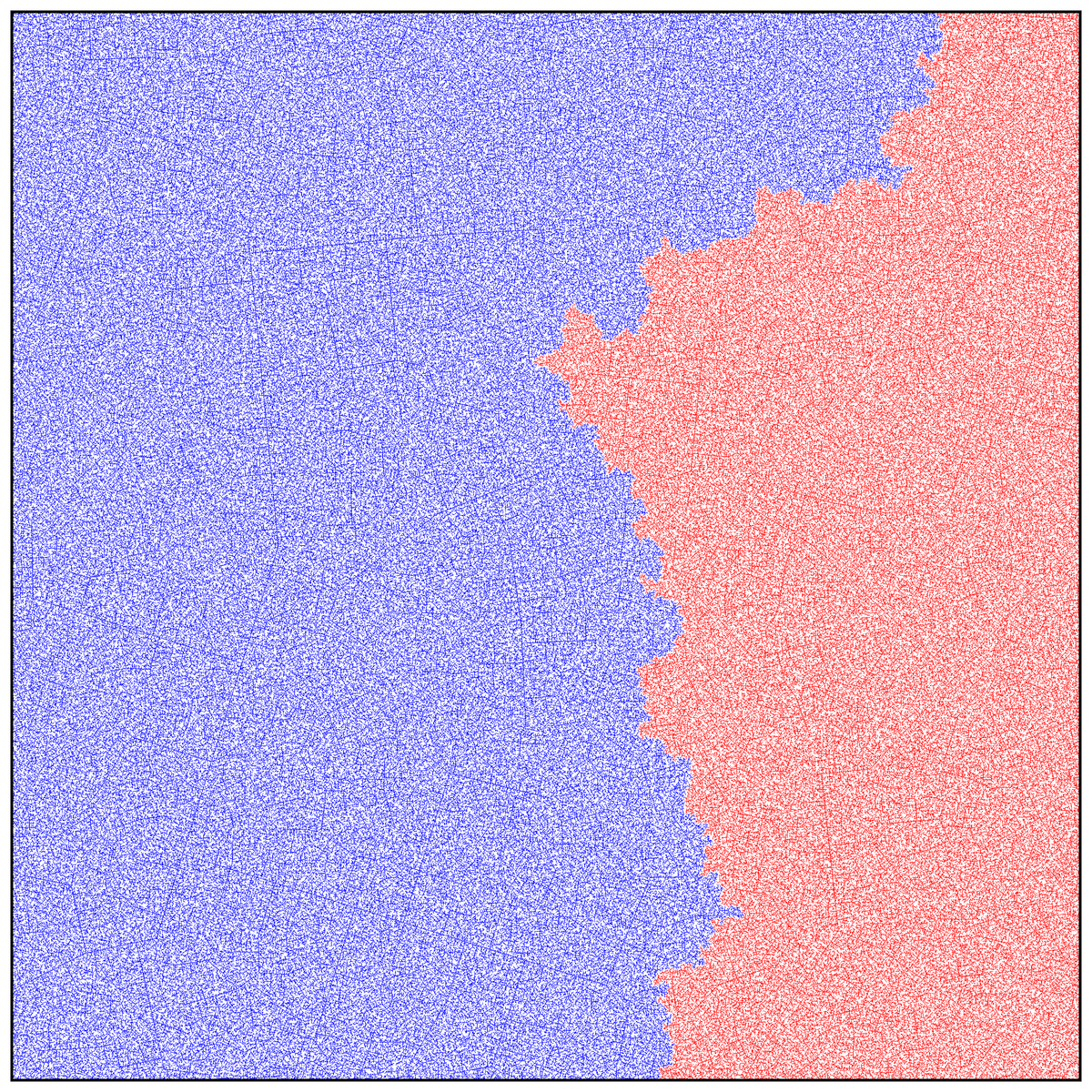}
\includegraphics[height=3cm]{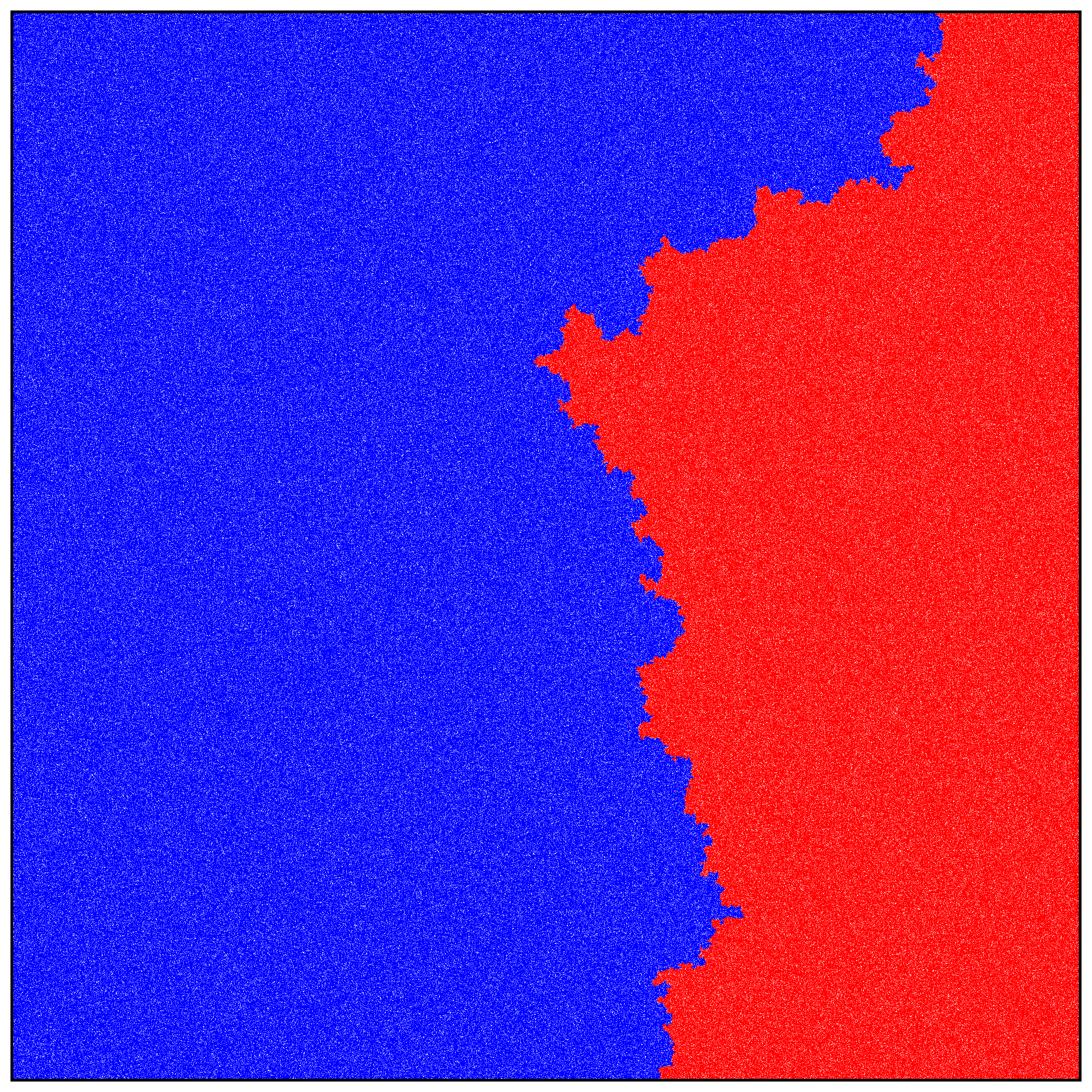}
\caption{Simulation of the variant Poisson coloring of space where a new incoming point is linked by a monochromatic segment to the nearest point in the process so far, from left to right with $10^{2}, 10^{3},10^{4},10^{6}$ and $10^{7}$ points. The arrivals points $X_i$ are the same as those used for Figure \ref{fig:simpoint}. Notice that the red ``island'' on the top right part of the figure present in the original point model has disappeared.}
\end{center}
\end{figure}

\medskip

\begin{figure}[ht]
\begin{center}
\includegraphics[width=15cm]{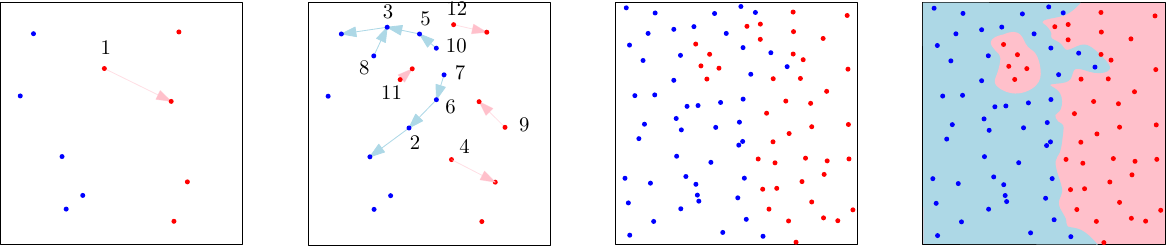}
\caption{Illustration of the creation of an ``island''. Such an island can be seen on Figure \ref{fig:simpoint} on the top right corner. This shows that the limiting frontier $ \mathcal{F}_\infty$ is not connected almost surely. \label{fig:island}}
\end{center}
\end{figure}

\textbf{Elusive topological properties of the frontier.}  
Although our results show the convergence in a strong sense of the colored regions and establish the fractal nature of the frontier, many questions remain open, such as the existence of a zero--one law for the Hausdorff dimension of $ \mathcal{F}_\infty$. 
We focus here  on  the planar case $d=2$, which concentrates the most interesting topological questions. 
Notice first that almost surely, the frontier $ \mathcal{F}_\infty$ is not connected, the reason being that it is possible for a point to get surrounded by points of the opposite color, thus eventually creating an ``island'' in the coloring. 
See \cite[Theorem 3]{preater2009species}, and Figure \ref{fig:island} below for an illustration. 
This island creation is not possible in the segment model, where the limiting frontier is almost surely connected. 

\begin{open}[Curves] Is the frontier $ \mathcal{F}_\infty$ a countable union of curves? \label{open:curve} \end{open}

It is natural to believe that $ \mathcal{F}_\infty$ is a countable union of curves \big(i.e, images of continuous paths from $[0,1]$ to $\mathbb{R}^2$\big), or that the limiting frontier in the segment model is a curve.
Although \cite{aizenman1999holder} provide sufficient conditions (namely \cite[Hypothesis \textbf{H1}]{aizenman1999holder}) which would allow to show that $\mathcal{F}_\infty$ \emph{contains} a curve\footnote{There are connected compact subsets of $ \mathbb{R}^2$ which do not contain any non-trivial curves, such as the pseudo-arc.}, checking those estimates seems hard in our setup due to the lack of a correlation inequality.
Yet, simulations suggest that the connected components of $ \mathcal{F}_\infty$ are \emph{simple} curves, meaning that ``double points'', i.e, points from which four alternating monochromatic non-trivial curves originate, do not exist.
\begin{open}[Simple curves] If the above question has a positive answer, are those curves almost surely \emph{simple}?\end{open}
\begin{figure}
\begin{center}
\includegraphics[height=4cm]{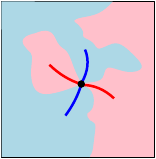}  \hspace{1cm}
\includegraphics[height=4cm]{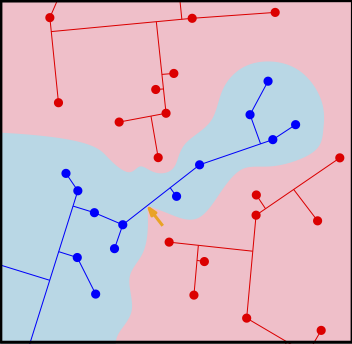}
\caption{\textbf{Left.} Illustration of a double point in the frontier $ \mathcal{F}_{\infty}$. \textbf{Right.} Can the frontier intersect the finite trees in the segment model (orange arrow)?}
\end{center}
\end{figure}

Accordingly, if that is true, then the frontier in the segment model should be made of a single simple curve. 
In fact, simulations suggest that in this model, the finite red and blue trees $ \mathcal{R}_n$ and $  \mathcal{B}_n$ are in the interior of the limiting red and blue regions $ \mathcal{R}_\infty$ and $ \mathcal{B}_\infty$ (it is possible to show that the arrival vertices $X_1,X_2,\ldots$ are indeed in the interior of $ \mathcal{R}_\infty$ and $ \mathcal{B}_\infty$, with minor technical modifications in the proof of Lemma \ref{lem:isolated} below, but the same results for the whole segments is still out of scope.) A more general question is the following.

\begin{open}[Safe margin] 
Suppose that $ \mathcal{R}_0$ is made of a segment or a ball instead of a single point. 
Do we have $ \mathbb{P}( \mathcal{R}_0  \cap \mathcal{B}_\infty\neq\emptyset) = 0$?
\end{open}

Our techniques (or those of Preater) only show that the above probability is strictly smaller than $1$, see the discussion before Corollary 1 in \cite{preater2009species}.

\section{Monochromatic balls and upper bound on $\dim_H\mathcal{F}_\infty$}

In this section, we establish our key lemma (Lemma \ref{lem:isolated}), which shows that for every $x\in[0,1]^d$ and $r>0$ such that $\overline{B}(x,r)$ does not contain the seeds $R_0$ and $B_0$, there is a positive probability that the smaller ball $\overline{B}(x,r/6)$ is monochromatic at the end of the coloring.
Applying Lemma \ref{lem:isolated} at all scales yields the upper bound on the dimension of the frontier.
In particular, it shows that for every $x\in[0,1]^d$, almost surely there exists an $r>0$ such that the ball $\overline{B}(x,r)$ is monochromatic at the end of the coloring.

\subsection{Key lemma}

Before stating the result, let us embed the model in continuous time to gain convenient independence properties.

\medskip

\textbf{Poissonization.}
Let $\Lambda$ be a Poisson process with intensity $\lambda\otimes\lambda_d$ on $\mathbb{R}_+\times\mathbb{R}^d$, where $\lambda$ and $\lambda_d$ denote the Lebesgue measures on $\mathbb{R}_+$ and $\mathbb{R}^d$, respectively. Let $X_1,X_2,\ldots$ be the points of $\Lambda$ that fall successively in $[0,1]^d$, at times say $\tau_1 < \tau_2 <\ldots$.
It is a standard fact that the $(X_n)_{n\in\mathbb{N}^*}$ are independent random variables, uniformly distributed in $[0,1]^d$.
Now, the coloring process can be defined in continuous time as follows.
The sequence $\left(\mathcal{R}_n,\mathcal{B}_n\right)_{n\in\mathbb{N}^*}$ of the discrete setting will here correspond to $\left(\mathcal{R}_{\tau_n},\mathcal{B}_{\tau_n}\right)_{n\in\mathbb{N}^*}$, and the sets $\mathcal{R}_t$ and $\mathcal{B}_t$ will be defined at all times $t\in\mathbb{R}_+$ as follows: for each $n\in\mathbb{N}$, let $\mathcal{R}_t=\mathcal{R}_{\tau_n}$ and $\mathcal{B}_t=\mathcal{B}_{\tau_n}$ for all $t\in[\tau_n,\tau_{n+1}[$, with the convention $\tau_0=0$.

\medskip

\begin{figure}
\begin{center}
\includegraphics[width=16cm]{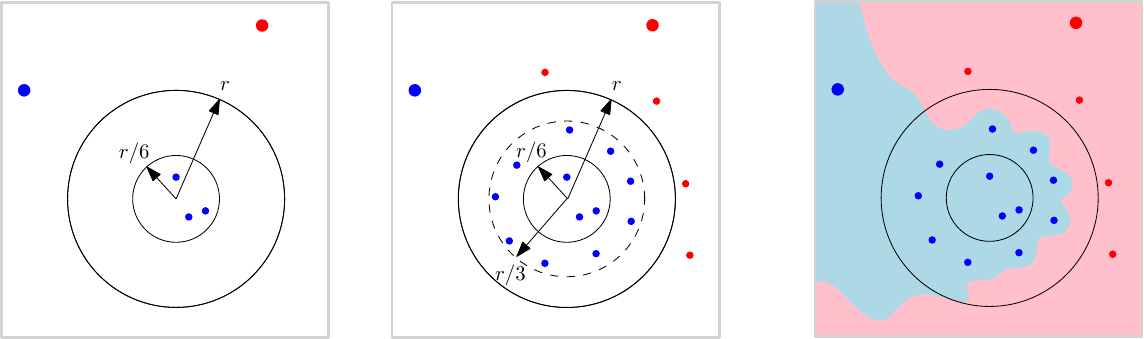}
\caption{Illustration of the construction of the event $G^x_{r/6,r}$.  Suppose that the seeds $R_{0}, B_{0}$ lie outside $\overline{B}(x,r)$. \\
\textbf{Left.} Notice that before the first point falls in $\overline{A}(x;r/6,r)$, a few points may have fallen in $\overline{B}(x,r/6)$ but, by the choice of $r/6$ versus $r$, all these points will have the same color (blue here).\\
\textbf{Center.} Then, we ask that the first few points falling in $\overline{A}(x;r/6,r)$ fall inside $\overline{B}(x,r/3)$, so that they all take the same color (blue here), and are spread out well enough to protect $\overline{B}(x,r/6)$ from being invaded by points of the other color (red here). \\\textbf{Right.} A reinforcement property of the process will then entail that with positive probability, the invaders (the red here) cannot penetrate $\overline{B}(x,r/6)$, so that it remains monochromatic. \label{fig:defense}}
\end{center}
\end{figure}

For $x\in\mathbb{R}^d$ and $0<r<R$, we denote by $\overline{A}(x;r,R)$ the annulus $\left\{y\in\mathbb{R}^d:r<|x-y|\leq R\right\}$, and we let $\mathcal{A}^x_{r,R}$ be the $\sigma$-algebra generated by the restriction of $\Lambda$ to $\mathbb{R}_+\times\overline{A}(x;r,R)$. The point of Lemma \ref{lem:isolated} below is to describe an $\mathcal{A}^x_{r/6,r}$-measurable ``good event'' $G^x_{r/6,r}$, which has probability bounded away from $0$ uniformly in $x$ and $r$, such that if $\overline{B}(x,r)$ does not contain the seeds $R_0$ and $B_0$, then on $G^x_{r/6,r}$ the ball $\overline{B}(x,r/6)$ is monochromatic at the end of the coloring. Figure \ref{fig:defense} provides an overview of how such a good event is constructed. 

\begin{lemma} \label{lem:isolated}
There exists a constant $p\in{]0,1[}$ for which the following holds.
For every $x\in[0,1]^d$ and $r>0$, there exists an $\mathcal{A}^x_{r/6,r}$-measurable good event $G^x_{r/6,r}$, which has probability $\mathbb{P}\left(G^x_{r/6,r}\right)\geq p$, such that if $\overline{B}(x,r)$ does not contain either $R_0$ or $B_0$, then on $G^x_{r/6,r}$ the ball $\overline{B}(x,r/6)$ does not meet both $\bigcup_{t\geq0}\mathcal{R}_t$ and $\bigcup_{t\geq0}\mathcal{B}_t$.
\end{lemma}
\begin{remark}\label{remkeylemma}
It will be clear from the proof that the event $G^x_{r/6,r}$ also prevents the ball $\overline{B}(x,r/6)$ from bichromaticity whenever $R_0\in\overline{B}(x,r/3)$ and $B_0\notin\overline{B}(x,r)$, or $B_0\in\overline{B}(x,r/3)$ and $R_0\notin\overline{B}(x,r)$. In particular, Lemma \ref{lem:isolated} allows to recover the result of \cite[proof of Theorem 2]{preater2009species} that almost surely, there exists an $r>0$ such that $\overline{B}(R_0,r)$ does not contain a blue point, and $\overline{B}(B_0,r)$ does not contain a red point.
\end{remark}
\begin{proof}
Fix $x\in[0,1]^d$ and $r>0$, and suppose that both $R_0$ and $B_0$ lie outside $\overline{B}(x,r)$.
We construct an $\mathcal{A}^x_{r/6,r}$-measurable good event $G$ on which a ``defense'' is organized inside the annulus $\overline{A}(x;r/6,r)$, preventing $\overline{B}(x,r/6)$ from meeting both $\bigcup_{t\geq0}\mathcal{R}_t$ and $\bigcup_{t\geq0}\mathcal{B}_t$.

\medskip

\textsc{Definition of $G$}\textbf{.}
Let $\rho_k=\left(1+2^{-k}\right)\cdot r/6$ for all $k\in\mathbb{N}$, and let $(t_k)_{k\in\mathbb{N}}$ be a sequence of positive real numbers to be adjusted later, with ${T_k:=t_0+\ldots+t_k\rightarrow\infty}$ as $k\to\infty$.
We define $\mathcal{A}^x_{r/6,r}$-measurable events $(G_k)_{k\in\mathbb{N}}$ such that for every $k\in\mathbb{N}$, on ${G_0\cap\ldots\cap G_k}$ the ball $\overline{B}(x,\rho_k)$ does not meet both $\bigcup_{0\leq t<T_k}\mathcal{R}_t$ and $\bigcup_{0\leq t<T_k}\mathcal{B}_t$.
The good event $G$ will then be defined as ${G=\bigcap_{k\geq0}G_k}$.
For every ${k\in\mathbb{N}}$, we denote by $A_k$ the annulus $\overline{A}(x;\rho_{k+1},\rho_k)$.
Let ${\delta_k=1/2\cdot(\rho_k-\rho_{k+1})\cdot(k+1)^{-2}}$, and let $\mathcal{Z}_k\subset A_k$ be a finite set of points with the following properties:
\begin{enumerate}[label=\alph*.]
\item\label{pZ1} for every $y\in A_k$, there exists $z\in\mathcal{Z}_k$ such that $y\in\overline{B}(z,3\delta_k/2)$,
\item\label{pZ2} for any $z\neq z'\in\mathcal{Z}_k$, we have $|z-z'|>\delta_k$,
\item\label{pZ3} for every $z\in\mathcal{Z}_k$, we have $\overline{B}(z,\delta_k/2)\subset A_k$.
\end{enumerate}
It is clear that such a set $\mathcal{Z}_k$ always exists: we keep adding points satisfying \ref{pZ2} and \ref{pZ3} until no more point can be added and then \ref{pZ1} must also be satisfied by construction. 
Note also that, because the balls $\left(\overline{B}(z,\delta_k/2)\right)_{z\in\mathcal{Z}_k}$ are disjoint and included in $A_k\subset\overline{B}(x,r/3)$, we have (by a volume argument)
\begin{equation}\label{Zkbound}
\#\mathcal{Z}_k\leq\left(\frac{r/3}{\delta_k/2}\right)^d=\left(8\cdot 2^{k+1}\cdot(k+1)^2\right)^d.
\end{equation}

We define $G_0$ as the event: ``for every $z\in\mathcal{Z}_0$, a point of $\Lambda$ falls in $\overline{B}(z,\delta_0/2)$ over the time interval $[0,t_0[$, meanwhile no point falls in $\overline{A}(x;r/3,r)$''.
We claim that on $G_0$, the ball $\overline{B}(x,r/3)$ does not meet both $\bigcup_{0\leq t<t_0}\mathcal{R}_t$ and $\bigcup_{0\leq t<t_0}\mathcal{B}_t$.
In particular, all the points of $\Lambda$ that have fallen in the spots $\left(\overline{B}(z,\delta_0/2)\right)_{z\in\mathcal{Z}_0}$ over the time interval $[0,t_0[$ have the same \emph{good} color.
Indeed, fix a realization of the event $G_0$.
Denote by $y_1,\ldots,y_n$ the points of $\Lambda$ that fall in $\overline{B}(x,r)$ over the time interval $[0,t_0[$, and by $\tau_1<\ldots<\tau_n\in[0,t_0[$ their arrival times.
Note that by the definition of $G_0$, the points $y_1,\ldots,y_n$ land in $\overline{B}(x,r/3)$.
So $y_1$ arrives in $\overline{B}(x,r/3)$, with its color.
Then when $y_2$ arrives, it lands within distance $2r/3$ of $y_1$, and at distance more than $2r/3$ of any other point of the process, since these all lie outside $\overline{B}(x,r)$.
Therefore, the nearest neighbor of $y_2$ is $y_1$, and $y_2$ inherits it color.
The argument iterates, proving the claim.

\begin{figure}[ht]
\begin{center}
\begin{tabular}{lr}
\includegraphics[width=0.4\linewidth]{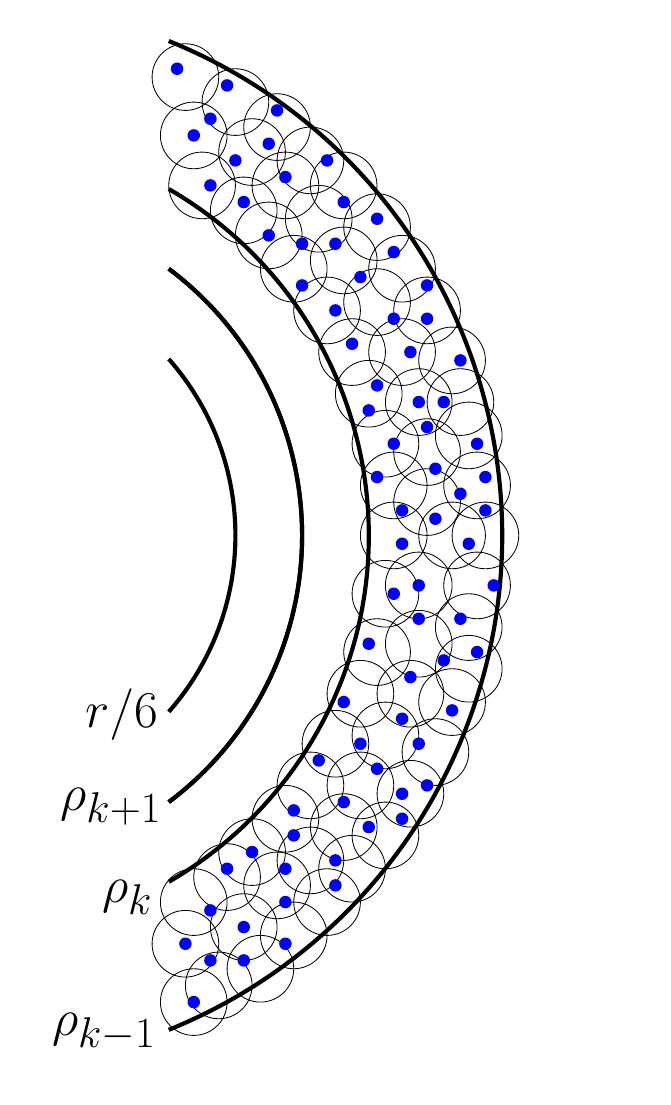}&\includegraphics[width=0.4\linewidth]{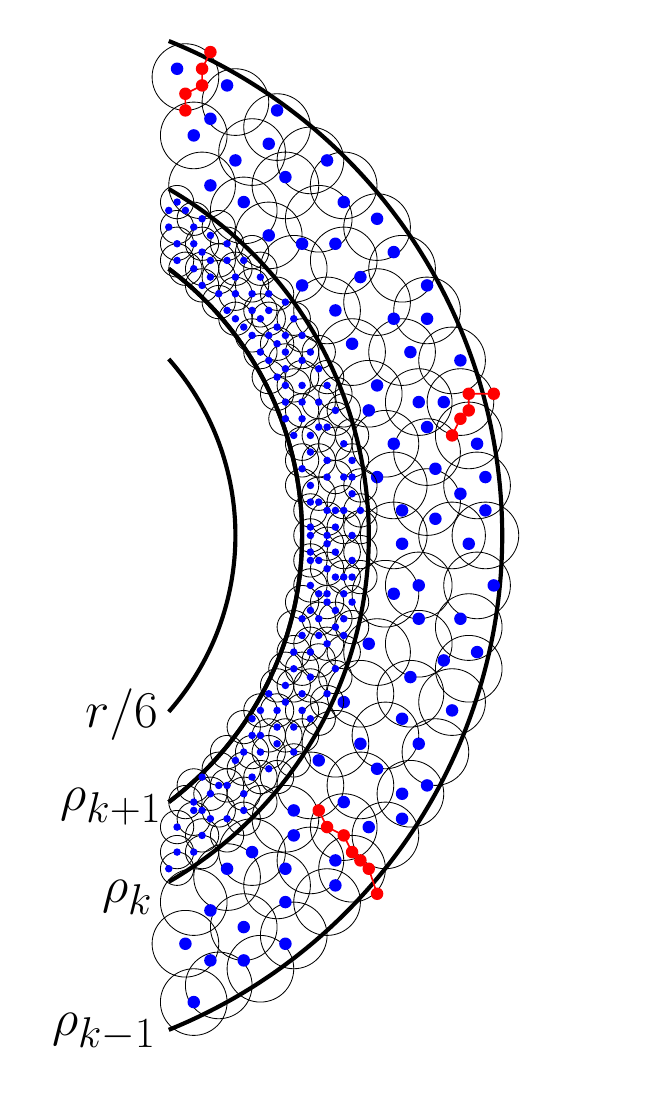}
\end{tabular}
\caption{Schematic description of the induction procedure. \\
\textbf{Left.} Assume that each cell of size $\delta_{k-1}$ in the annulus between radii $\rho_{k-1}$ and $\rho_k$ around $x$ contains a blue point, and that no red point has penetrated the ball with radius $\rho_{k-1}$ around $x$.\\
\textbf{Right.} Then, if red points want to cross the previous annulus, they need to navigate between the blue points present in that annulus, which forces them to create a chain of length at least $k^2$. Although it may look the other way around on the picture (due to the difficulty of representing realistic scales), the following is true: by the time such a chain is created, with very high probability, each cell of size $\delta_k$ in the annulus between radii $\rho_k$ and $\rho_{k+1}$ around $x$ will contain a blue point.}
\end{center}
\end{figure}

Next, in order to define $G_k$ for $k\geq1$, we start with the following deterministic observation.
Suppose by induction that, on the event $G_0\cap\ldots\cap G_{k-1}$, the following holds: at time $T_{k-1}$, each cell $\left(\overline{B}(z,\delta_{k-1}/2)\right)_{z\in\mathcal{Z}_{k-1}}$ contains a point of the good color, and $\overline{B}(x,\rho_{k-1})$ does not contain any point of the other bad color.
Then, every $y\in A_{k-1}$ is within distance $2\delta_{k-1}$ of a point of the good color, and the only way of bringing a point of the bad color inside $\overline{B}(x,\rho_k)$ before time $T_k$ is to have points, say $y_1,\ldots,y_j$, of $\Lambda$ falling in $A_{k-1}$, at times say $\tau_1<\ldots<\tau_j\in[T_{k-1},T_k[$, with:
\begin{itemize}
\item $d\left(y_1;\left.\mathbb{R}^d\middle\backslash\overline{B}(x,\rho_{k-1})\right.\right)<2\delta_{k-1}$,
\item $|y_{i+1}-y_i|<2\delta_{k-1}$ for each $i\in\llbracket1,j\llbracket$,
\item $d\left(y_j,\overline{B}(x,\rho_k)\right)<2\delta_{k-1}$.
\end{itemize}
Now, let us discretise this information.
First, it follows from the inequality $\rho_{k-1}-\rho_k<(j+1)\cdot2\delta_{k-1}$ that such a path must have length $j\geq k^2$. Then, for each $i\in\left\llbracket1,k^2\right\rrbracket$, let $z_i\in\mathcal{Z}_{k-1}$ be such that ${y_i\in\overline{B}(z_i,3\delta_{k-1}/2)}$.
The following holds:
\begin{itemize}
\item for every $i\in\left\llbracket1,k^2\right\rrbracket$, we have $z_i\in\mathcal{Z}_{k-1}$,
\item for each $i\in\left\llbracket1,k^2\right\llbracket$, we have $|z_{i+1}-z_i|\leq5\delta_{k-1}$.
\end{itemize}
A sequence $z_1,\ldots,z_{k^2}$ satisfying the two properties above is said to be \emph{admissible of order $k$}. 
Moreover, since for each $i\in\left\llbracket1,k^2\right\rrbracket$, a point of $\Lambda$ falls in $\overline{B}(z_i,3\delta_{k-1}/2)$ at time $\tau_i$, with 
\[\tau_1<\ldots<\tau_{k^2}\in[T_{k-1},T_k[,\]
we say that $z_1,\ldots,z_{k^2}$ \emph{ring consecutively} over the time interval ${[T_{k-1},T_k[}$.
We can now formally define the event $G_k$ by ``for every $z\in\mathcal{Z}_k$, a point of $\Lambda$ falls in $\overline{B}(z,\delta_k/2)$ over the time interval $[T_{k-1},T_k[$, meanwhile no admissible sequence of order $k$ rings consecutively''. By induction on $k$, we see that on $G_0\cap\ldots\cap G_k$, the ball $\overline{B}(x,\rho_k)$ does not meet both $\bigcup_{0\leq t<T_k}\mathcal{R}_t$ and $\bigcup_{0\leq t<T_k}\mathcal{B}_t$. Finally, we set $G=\bigcap_{k\geq0}G_k$.

\medskip

\textsc{The probability of $G$ is bounded away from $0$}\textbf{.}
Because of the disjointness of the time intervals over which they are defined, the events $(G_k)_{k\in\mathbb{N}}$ are independent:
\[\mathbb{P}(G)=\prod_{k\geq0}\mathbb{P}(G_k)=\mathbb{P}(G_0)\cdot\prod_{k\geq1}[1-\mathbb{P}(F_k)],\]
where $F_k$ is the complement of the event $G_k$.
Now, on $F_k$, we have the following alternative.
\begin{itemize}
\item Either there exists $z\in\mathcal{Z}_k$ such that $\Lambda\left({[T_{k-1},T_k[}\times\overline{B}(z,\delta_k/2)\right)=0$: let us call $B_k$ the corresponding event.
\item Or there exists an admissible sequence of order $k$ that rings consecutively: we call $C_k$ the corresponding event.
\end{itemize}
We have $\mathbb{P}(F_k)\leq\mathbb{P}(B_k)+\mathbb{P}(C_k)$.
For the second term, the union bound and the Markov property for $\Lambda$ show that
\[\begin{split}
\mathbb{P}(C_k)&\leq\sum_\text{$z_1,\ldots,z_{k^2}$ admissible of order $k$}\mathbb{P}(\text{$z_1,\ldots,z_{k^2}$ ring consecutively})\\
&=\#\{\text{admissible sequences of order $k$}\}\cdot\mathbb{P}(\tau_1+\ldots+\tau_{k^2}<t_k),
\end{split}\]
where the $\tau_i$'s are independent exponential random variables with inverse mean ${\lambda_k=\upsilon_d\cdot(3\delta_{k-1}/2)^d}$, and $\upsilon_d$ is the Lebesgue measure of the unit ball in $\mathbb{R}^d$. 
We set $t_k=\alpha\cdot k^2\cdot\lambda_k^{-1}$ for all $k\geq1$, where $\alpha\in{]0,1[}$ is a parameter to be adjusted later. 
On the one hand, the Chernoff bound yields
\[\mathbb{P}(\tau_1+\ldots+\tau_{k^2}<t_k)\leq e^{(1-\alpha+\ln\alpha)k^2}.\]
On the other hand, in order to choose an admissible sequence $z_1,\ldots,z_{k^2}$ of order $k$, there is no more than $\#\mathcal{Z}_{k-1}$ possibilities for the choice of $z_1$, and then for each ${i\in\left\llbracket1,k^2\right\llbracket}$, there is at most $\#\left(\mathcal{Z}_{k-1}\cap\overline{B}(z_i,5\delta_{k-1})\right)\leq11^d$ possibilities for the choice of $z_{i+1}$; this last upper bound holds by a volume argument, since the balls $\left(\overline{B}(z,\delta_{k-1}/2)\,;\,z\in\mathcal{Z}_{k-1}\cap\overline{B}(z_i,5\delta_{k-1})\right)$ are disjoint and included in $\overline{B}(z_i,11\delta_{k-1}/2)$. 
Thus, we find that
\[\mathbb{P}(C_k)\leq\#\mathcal{Z}_{k-1}\cdot\left(11^d\right)^{k^2-1}\cdot e^{(1-\alpha+\ln\alpha)k^2}.\]
We now fix $\alpha$ such that $11^d\cdot e^{1-\alpha+\ln\alpha}\leq e^{-1}$. Recalling \eqref{Zkbound}, we obtain
\[\mathbb{P}(C_k)\leq\left(8\cdot 2^k\cdot k^2\right)^d\cdot11^{-d}\cdot e^{-k^2}=:c_k.\]
Next, we upper bound $\mathbb{P}(B_k)$. We have $\mathbb{P}(B_k)\le \#\mathcal{Z}_k \cdot p_k$, where
\[p_k=\exp\left[-\upsilon_d\cdot(\delta_k/2)^d\cdot t_k\right]=\exp\left[-\alpha\cdot\left(\frac{\delta_k}{3\delta_{k-1}}\right)^d\cdot k^2\right].\]
Using again \eqref{Zkbound}, we see that
\[\#\mathcal{Z}_k\cdot p_k\leq\left(8\cdot2^{k+1}\cdot(k+1)^2\right)^d\cdot\exp\left[-\alpha\cdot24^{-d}\cdot k^2\right]=:b_k,\]
which finally yields
\[\mathbb{P}(F_k)\leq b_k+c_k=:a_k.\]

Now, we can complete the proof that $\mathbb{P}(G)$ is bounded away from $0$ uniformly in $x$ and $r$. Since $\sum_{k\geq1}a_k<\infty$, we can find $K\in\mathbb{N}$ (not depending on $x$ or $r$) such that $\prod_{k\geq K+1}(1-a_k)\geq1/2$. With that choice, we have
\[\mathbb{P}(G)=\mathbb{P}(G_0)\cdot\ldots\cdot\mathbb{P}(G_K)\cdot\prod_{k\geq K+1}[1-\mathbb{P}(F_k)]\geq\frac{\mathbb{P}(G_0)\cdot\ldots\cdot\mathbb{P}(G_K)}{2}.\]
Note that we have yet to specify the value of $t_0$, which we now set to $t_0=r^{-d}$. Given this choice, the probability $\mathbb{P}(G_0)$ is bounded away from $0$ uniformly in $x$ and $r$. Next, for each $k\in\llbracket1,K\rrbracket$, we claim that $\mathbb{P}(G_k)$ is also bounded away from $0$ uniformly in $x$ and $r$, because the same is true for the probability of the sub-event: ``for each $z\in\mathcal{Z}_k$, a point of $\Lambda$ falls in $\overline{B}(z,\delta_k/2)$ over the time interval $[T_{k-1},T_k[$, meanwhile no point falls in $A_{k-1}$''. Thus, the quantity $\mathbb{P}(G_0)\cdot\ldots\cdot\mathbb{P}(G_K)\cdot1/2$ is bounded away from zero uniformly in $x$ and $r$, which concludes the proof of the lemma.
\end{proof}

\subsection{Hausdorff dimension: upper bound}

\begin{proposition}\label{propkeylemma'}
There exists constants $C, \alpha >0$ such that for every $x\in[0,1]^d$,
\[\mathbb{P}\left(\text{$\mathcal{F}_\infty$ meets $\overline{B}(x,\delta)$}\right)\leq C\cdot\delta^\alpha\quad\text{for all $\delta\in{]0,1[}$ small enough so that $R_0,B_0\notin\overline{B}\left(x,\sqrt{\delta}\right)$}.\]
\end{proposition}
\begin{proof}
Fix $x\in[0,1]^d$, and let $\delta\in{]0,1[}$ be small enough so that $R_0,B_0\notin\overline{B}\left(x,\sqrt{\delta}\right)$.
Set ${r_k=\sqrt{\delta}\cdot6^{-k}}$ for all $k\in\mathbb{N}$, and let $K=\left\lceil\log_6\left(\delta^{-1/2}\right)\right\rceil-1$ be the largest integer $k\in\mathbb{N}$ such that $r_k>\delta$.
By Lemma \ref{lem:isolated}, we have the inclusion
\[\left(\text{$\mathcal{F}_\infty$ meets $\overline{B}(x,\delta)$}\right)\subset\left(\text{for every $k\in\llbracket1,K\rrbracket$, the event $G^x_{r_k,r_{k-1}}$ fails to be realized}\right).\]
Thus, since those are independent events, we obtain
\[\mathbb{P}\left(\text{$\mathcal{F}_\infty$ meets $\overline{B}(x,\delta)$}\right)\leq(1-p)^K\leq(1-p)^{-1}\cdot(1-p)^\alpha,\]
where $\alpha=-\ln(1-p)/(2\ln6)>0$.
This is the desired upper bound.
\end{proof}

\begin{proposition}\label{prop:upperbound}
There exists $\varepsilon>0$ such that, almost surely,
\[\dim_H\mathcal{F}_\infty\leq d-\varepsilon<d.\]
\end{proposition}
\begin{proof}
Let $\varepsilon=\alpha\wedge(d/2)$, where $\alpha$ is the exponent of Proposition \ref{propkeylemma'}.
For each $k\in\mathbb{N}$, set $\delta_k=2^{-k}$, and let $\left(\overline{B}(x,\delta_k)\right)_{x\in\mathcal{X}_k}$ be a covering of $[0,1]^d$ by balls of radius $\delta_k$, with centers $x\in[0,1]^d$ more than $\delta_k$ apart so that the $\left(\overline{B}(x,\delta_k/2)\right)_{x\in\mathcal{X}_k}$ are disjoint.
In particular, there exists a constant ${C'=C'(d)>0}$ such that $\#\mathcal{X}_k\leq C'\cdot\delta_k^{-d}$.
By definition, the $(d-\varepsilon)$-dimensional Hausdorff measure of $\mathcal{F}_\infty$ is bounded from above by
\[H=\varliminf_{k\to\infty}\sum_{x\in\mathcal{X}_k}(2\delta_k)^{d-\varepsilon}\cdot\mathbf{1}\left(\text{$\mathcal{F}_\infty$ meets $\overline{B}(x,\delta_k)$}\right).\]
Now, we claim that the random variable $H$ is finite a.s, which implies that ${\dim_H\mathcal{F}_\infty\leq d-\varepsilon}$ almost surely. Indeed, using Fatou's lemma, we get
\[\mathbb{E}[H]\leq\varliminf_{k\to\infty}\sum_{x\in\mathcal{X}_k}(2\delta_k)^{d-\varepsilon}\cdot\mathbb{P}\left(\text{$\mathcal{F}_\infty$ meets $\overline{B}(x,\delta_k)$}\right),\]
and for each $k$ we have, using Proposition \ref{propkeylemma'}:
\begin{multline*}
\sum_{x\in\mathcal{X}_k}(2\delta_k)^{d-\varepsilon}\cdot\mathbb{P}\left(\text{$\mathcal{F}_\infty$ meets $\overline{B}(x,\delta_k)$}\right)\\
\leq\#\left\{x\in\mathcal{X}_k:\text{$R_0\in\overline{B}\left(x,\sqrt{\delta_k}\right)$ or $B_0\in\overline{B}\left(x,\sqrt{\delta_k}\right)$}\right\}\cdot(2\delta_k)^{d-\varepsilon}+\#\mathcal{X}_k\cdot(2\delta_k)^{d-\varepsilon}\cdot C\cdot\delta_k^\alpha.
\end{multline*}
For the first term, we have
\begin{eqnarray*}
\lefteqn{\#\left\{x\in\mathcal{X}_k:\text{$R_0\in\overline{B}\left(x,\sqrt{\delta_k}\right)$ or $B_0\in\overline{B}\left(x,\sqrt{\delta_k}\right)$}\right\}}\\
&\leq&\#\mathcal{X}_k\cap\overline{B}\left(R_0,\sqrt{\delta_k}\right)+\#\mathcal{X}_k\cap\overline{B}\left(B_0,\sqrt{\delta_k}\right)\\
&\leq&\left(\frac{\sqrt{\delta_k}}{\delta_k/2}+1\right)^d+\left(\frac{\sqrt{\delta_k}}{\delta_k/2}+1\right)^d=2\cdot\left(\frac{2}{\sqrt{\delta_k}}+1\right)^d.
\end{eqnarray*}
Recalling that $\varepsilon\leq d/2$, we deduce that
\[\varliminf_{k\to\infty}\#\left\{x\in\mathcal{X}_k:\text{$R_0\in\overline{B}\left(x,\sqrt{\delta_k}\right)$ or $B_0\in\overline{B}\left(x,\sqrt{\delta_k}\right)$}\right\}\cdot(2\delta_k)^{d-\varepsilon}<\infty.\]
For the second term, using that $\#\mathcal{X}_k\leq C'\cdot\delta_k^{-d}$ and $\varepsilon\leq\alpha$, we obtain that
\[\varliminf_{k\to\infty}\#\mathcal{X}_k\cdot(2\delta_k)^{d-\varepsilon}\cdot C\cdot\delta_k^\alpha<\infty.\]
Combining these two inequalities, we conclude that $\mathbb{E}[H]<\infty$, which completes the proof.
\end{proof}

\section{Hausdorff dimension lower bounds}

In this section, we prove that the Hausdorff dimension of $\mathcal{F}_\infty$ is strictly greater than $d-1$.
A substantial part of this work consists in adaptating the lower bound of \cite[Theorem 1.3]{aizenman1999holder}. Indeed, as the knowledgeable reader has undoubtedly noticed, it is not possible to invoke this result directly because we do not know whether the frontier $\mathcal{F}_\infty$ contains non-trivial curves. So instead, we modify the proof of Aizenman and Burchard to obtain a general Hausdorff dimension lower bound result for connected random closed subsets  which satisfy Property $(\varnothing)$, with the following definition.

\begin{definition}[Property $(\varnothing)$]
Let $\mathcal{F}$ be a random closed subset of $[0,1]^d$.
We say that $\mathcal{F}$ satisfies Property $(\varnothing)$ if there exists a constant $\zeta>1$, and two constants $Q>0$ and $q\in{]0,1[}$ such that the following holds: for every collection $\left(\overline{B}(x_i,r_i)\,;\,i\in\llbracket1,n\rrbracket\right)$ of balls with centers $x_1,\ldots,x_n\in[0,1]^d$ such that the dilated balls $\left(\overline{B}(x_i,\zeta r_i)\,;\,i\in\llbracket1,n\rrbracket\right)$ are disjoint (we say that the balls are $\zeta$-separated), we have
\[\mathbb{P}\left(\text{for each $i\in\llbracket1,n\rrbracket$, the set $\mathcal{F}$ meets $\overline{B}(x_i,r_i)$}\right)\leq Q\cdot q^n.\]
\end{definition}

We point out that Property $(\varnothing)$ is very similar to \cite[Hypothesis \textbf{H2}]{aizenman1999holder}.

\begin{theorem}\label{thmAB'}
Let $\mathcal{F}$ be a random closed subset of $[0,1]^d$ that satisfies Property $(\varnothing)$.
There exists a constant $s>1$ such that almost surely, the following holds: for every non-trivial (i.e, not empty or reduced to a point) connected closed subset $\mathcal{C}\subset\mathcal{F}$, we have
\[\dim_H\mathcal{C}\geq s>1.\]
\end{theorem}

The empty set, or the points of a homogeneous Poisson process, are obvious examples of random closed subsets of $[0,1]^d$ satisfying Property $(\varnothing)$.
Both have Hausdorff dimension $0$, but all their connected components are trivial.
The above result says that, as soon as we request a random closed subset to have a non-trivial connected component (and thus, Hausdorff dimension at least $1$), then the fact that it satisfies Property $ (\varnothing)$ implies that it is ``delocalized'' in some sense, and entails that its Hausdorff dimension is, in fact, strictly greater than $1$.

\medskip

In Subsection \ref{sec:lowerboundF}, we obtain the lower bound of Theorem \ref{thm:main} by applying Theorem \ref{thmAB'} to the frontier $\mathcal{F}_\infty$.
This directly gives the result $\dim_H\mathcal{F}_\infty>1$ a.s.~in dimension $d=2$; and the lower bound $\dim_H\mathcal{F}_\infty>d-1$ a.s.~in general dimensions follows from Theorem \ref{thmAB'} together with a slicing lemma (see Proposition \ref{prop:lowerbound} below).
Let us now present the proof of Theorem \ref{thmAB'}.

\subsection{Proof of Theorem \ref{thmAB'}}

As mentioned before, the proof is adapted from \cite{aizenman1999holder}, and thus uses similar ingredients. Still, we provide here a self-contained proof, recalling and adapting the necessary results from \cite{aizenman1999holder} whenever required. At its core, the proof employs the usual ``energy method'' (see \cite[Theorem 6.4.6]{bishop_peres_2016}) to lower bound the Hausdorff dimension of a set. There are two main parts:
\begin{enumerate}
\item We first describe a deterministic splitting procedure for curves which produces, when the curves are oscillating enough, an important number of disjoint sub-curves.
\item Next, we show that if a connected random closed subset $\mathcal{F}$ satisfies Property $(\varnothing)$, then curves located in shrinking neighborhoods of $\mathcal{F}$ will necessarily oscillate enough so that we 
can use the splitting procedure above to create many sub-curves. This will enable us to create a sequence of measures with good integrability properties and finally, by compactness, extract a measure $\nu$ supported on $\mathcal{F}$  such that $\iint|x-y|^{-s}\mathrm{d}\nu(y)\mathrm{d}\nu(x)<\infty$ for some $s>1$, which in turn implies that $\dim_H\mathcal{F}\geq s$.
\end{enumerate}

\subsubsection{A deterministic splitting procedure for curves}

\begin{figure}
\begin{center}
\includegraphics[width=12cm]{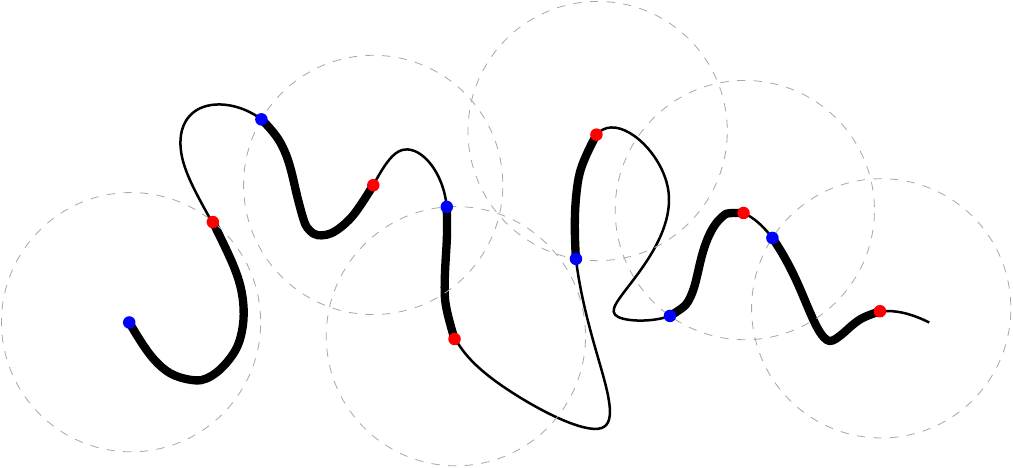}
\caption{Illustration of the procedure $( \mathrm{P}_{\alpha})$. The times $\sigma_{i}$ are in blue, and the $\tau_{i}$ in red. Several subpaths (in thick line), spanning a distance $\delta$ and being $\alpha \delta$ apart, are created from the initial path. The balls have radius $\delta$. \label{fig:procedure}}
\end{center}
\end{figure}

Given a small parameter $\alpha\in{]0,1[}$, let us describe the splitting procedure $(\mathrm{P}_\alpha)$ introduced in \cite[Lemma 5.2]{aizenman1999holder}. 
It takes as input a continuous path ${\gamma:[0,1]\rightarrow\mathbb{R}^d}$ with $\gamma(1)\neq\gamma(0)$, and outputs a collection $\gamma_1,\ldots,\gamma_\kappa$ of subpaths of $\gamma$, with the following properties:
\begin{itemize}
\item for every $i\in\llbracket1,\kappa\rrbracket$, we have $|\gamma_i(0)-\gamma_i(1)|=\alpha\cdot|\gamma(0)-\gamma(1)|=:\delta$,
\item for any $i\neq j\in\llbracket1,\kappa\rrbracket$, we have $d(\gamma_i[0,1];\gamma_j[0,1])\geq\alpha\delta$.
\end{itemize}
The splitting procedure $(\mathrm{P}_\alpha)$ goes as follows.
See Figure \ref{fig:procedure} for an illustration.
Let
\[\Delta=|\gamma(0)-\gamma(1)|>0,\quad\text{and}\quad\delta=\alpha\cdot\Delta<|\gamma(0)-\gamma(1)|.\]
Initially, set ${\sigma_1=0}$, and let 
\[\tau_1=\inf\left\{t\in[0,1]:\gamma(t)\notin\overline{B}(\gamma(0),\delta)\right\}.\]
Then, for $i\in\mathbb{N}^*$ such that $\sigma_1,\tau_1;\ldots;\sigma_i,\tau_i$ have been constructed: if
\[d(\gamma(t);\gamma[\sigma_1,\tau_1]\cup\ldots\cup\gamma[\sigma_i,\tau_i])\leq(1+\alpha)\delta\quad\text{for all $t\in[\tau_i,1]$},\]
then we set $\sigma_{i+1}=1$ and $\tau_{i+1}=1$. Otherwise, we set
\[\tau_{i+1}=\inf\{t\in[\tau_i,1]:d(\gamma(t);\gamma[\sigma_1,\tau_1]\cup\ldots\cup\gamma[\sigma_i,\tau_i])>(1+\alpha)\delta\},\]
and let $\sigma_{i+1}=\sup\left\{t\in{[\tau_i,\tau_{i+1}[}:\gamma(t)\notin\overline{B}(\gamma(\tau_{i+1}),\delta)\right\}$.
Finally, let $\kappa$ be the largest integer ${i\in\mathbb{N}^*}$ such that $\tau_i<1$, and for each $i\in\llbracket1,\kappa\rrbracket$, denote by $\gamma_i$ the path ${\theta\in[0,1]\mapsto\gamma((1-\theta)\sigma_i+\theta\tau_i)}$.
For every $i\in\llbracket1,\kappa\rrbracket$, we have $|\gamma_i(0)-\gamma_i(1)|=\delta$, and for any $i\neq j\in\llbracket1,\kappa\rrbracket$, we have ${d(\gamma_i[0,1];\gamma_j[0,1])\geq\alpha\delta}$.
 
\begin{definition}
We say that a continuous path $\gamma:[0,1]\rightarrow\mathbb{R}^d$, with $\gamma(1)\neq\gamma(0)$, \emph{deviates by a factor ${\rho>0}$ from being a straight line} when there exists $t\in[0,1]$ such that 
\[\gamma(t)\notin S(\gamma(0),\gamma(1);\rho|\gamma(0)-\gamma(1)|),\]
where $S(x,y;r)=\left\{z\in\mathbb{R}^d:d(z;[x,y])\leq r\right\}$ denotes the sausage of radius $r$ around the line segment $[x,y]$.
\end{definition}

Intuitively (recall that $\alpha$ is small) the number $\kappa$ of subpaths produced by the procedure $(\mathrm{P}_\alpha)$ must be at least of order ${\Delta/\delta=1/\alpha}$, and this lower bound is roughly attained by a straight line. 
However, when the input path deviates from being a straight line, one would expect the procedure to produce additional paths. 
This is the content of the next proposition.

\begin{proposition}\label{propkappa}
Let $\gamma:[0,1]\rightarrow\mathbb{R}^d$ be a path with $\gamma(1)\neq\gamma(0)$.
\begin{enumerate}
\item The number $\kappa$ of subpaths produced by the procedure $(\mathrm{P}_\alpha)$ always satisfies
\[\kappa\geq\frac{1-\alpha}{(1+\alpha)\alpha}.\]

\item If $\gamma$ deviates by a factor $\rho>0$ from being a straight line, then the number $\kappa$ of subpaths produced by the procedure $(\mathrm{P}_\alpha)$ satisfies
\[\kappa\geq\frac{\left.\left(1+\sqrt{1+4\rho^2}\right)\middle/2\right.-(4+\alpha)\alpha}{(1+\alpha)\alpha}.\]
\end{enumerate}
\end{proposition}

\begin{proof}
We keep the notation introduced above: we have $\Delta=|\gamma(0)-\gamma(1)|>0$, and  $\delta=\alpha\cdot\Delta$.
\begin{enumerate}
\item By the definition of $(\mathrm{P}_\alpha)$, there exists $i_1\in\llbracket1,\kappa\rrbracket$ such that ${d(\gamma(1);\gamma[\sigma_{i_1},\tau_{i_1}])\leq(1+\alpha)\delta}$, and $t_1\in[\sigma_{i_1},\tau_{i_1}]$ such that $|\gamma(1)-\gamma(t_1)|\leq(1+\alpha)\delta$.
Then, by induction, for $k\in\mathbb{N}^*$ such that $i_1,t_1;\ldots;i_k,t_k$ have been constructed, proceed as follows.
If $i_k=1$, then set $i_{k+1}=1$ and let $t_{k+1}=0$.
Otherwise, by the definition of $(\mathrm{P}_\alpha)$, there exists ${i_{k+1}\in\llbracket1,i_k\llbracket}$ such that $d\left(\gamma(t_k);\gamma\left[\sigma_{i_{k+1}},\tau_{i_{k+1}}\right]\right)\leq(1+\alpha)\delta$, and ${t_{k+1}\in\left[\sigma_{i_{k+1}},\tau_{i_{k+1}}\right]}$ such that 
\[|\gamma(t_k)-\gamma(t_{k+1})|\leq(1+\alpha)\delta.\]
Finally, let $m$ be the smallest integer $k\in\mathbb{N}^*$ such that $i_k=1$.
We have
\[\begin{split}
|\gamma(0)-\gamma(1)|&\leq|\gamma(0)-\gamma(t_m)|+\sum_{k=1}^{m-1}|\gamma(t_{k+1})-\gamma(t_k)|+|\gamma(t_1)-\gamma(1)|\\
&\leq\delta+m\cdot(1+\alpha)\delta,
\end{split}\]
hence
\[m\geq\frac{\Delta-\delta}{(1+\alpha)\delta}=\frac{1-\alpha}{(1+\alpha)\alpha}.\]
The result follows, since $i_1,\ldots,i_m$ are distinct elements of $\llbracket1,\kappa\rrbracket$.

\item Suppose that there exists $t\in[0,1]$ such that ${\gamma(t)\notin S(\gamma(0),\gamma(1);\rho\Delta)}$. 
We still denote by $i_1,t_1;\ldots;i_m,t_m$ the sequence of indices and times defined above. 
Moreover, we construct another sequence $j_1,u_1;\ldots;j_p,u_p$ in the exact same way but now obtained by backtracking from time $t$ instead of time $1$. 
By construction, we have  ${|\gamma(t)-\gamma(u_1)|\leq(1+\alpha)\delta}$, and ${|\gamma(u_k)-\gamma(u_{k+1})|\leq(1+\alpha)\delta}$ for all $k\in\llbracket1,p\llbracket$. 
Finally, let $n$ be the smallest integer $k\in\llbracket1,p\rrbracket$ such that $j_k\in \{i_1,\ldots,i_m\}$, and let $l\in\llbracket1,m\rrbracket$ be the unique index such that $j_n=i_l$.
By construction, the indices $i_1,\ldots,i_m$ and $j_1,\ldots,j_{n-1}$ are all distinct, hence $\kappa\geq m+n-1$.
Now, on the one hand, with the same argument as above, we have:
\[|\gamma(0)-\gamma(t_{l})|+|\gamma(t_{l})-\gamma(1)|\leq\delta+m\cdot(1+\alpha)\delta.\]
On the other hand, using that $\gamma(t_l)$ and $\gamma(u_n)$ both lie in $\gamma\left[\sigma_{i_l},\tau_{i_l}\right]$, which is contained in some ball of radius $\delta$ by construction, we have
\[\begin{split}
|\gamma(t_{l})-\gamma(t)|&\leq|\gamma(t_{l})-\gamma(u_n)|+\sum_{k=1}^{n-1}|\gamma(u_{k+1})-\gamma(u_k)|+|\gamma(u_1)-\gamma(t)|\\
&\leq2\delta+n\cdot(1+\alpha)\delta.
\end{split}\]
Summing these inequalities, we get
\[|\gamma(0)-\gamma(t_l)|+|\gamma(1)-\gamma(t_l)|+|\gamma(t)-\gamma(t_l)|\leq(m+n)(1+\alpha)\delta+3\delta,\]
hence
\[m+n-1\geq\frac{\inf_{x\in\mathbb{R}^d}\{|\gamma(0)-x|+|\gamma(1)-x|+|\gamma(t)-x|\}-(4+\alpha)\delta}{(1+\alpha)\delta}.\]
It remains to lower bound the infimum in the right hand side. First, using the triangle inequality, we get
\[|\gamma(0)-x|+|\gamma(1)-x|+|\gamma(t)-x|\geq\frac{|\gamma(0)-\gamma(1)|+|\gamma(0)-\gamma(t)|+|\gamma(t)-\gamma(1)|}{2}.\]
Then, we make use of the fact that $\gamma(t)\notin S(\gamma(0),\gamma(1);\rho\Delta)$, to get
\[|\gamma(0)-\gamma(t)|+|\gamma(t)-\gamma(1)|\geq\sqrt{1+4\rho^2}\cdot\Delta.\]
Altogether, we obtain
\[\inf_{x\in\mathbb{R}^2}\{|\gamma(0)-x|+|\gamma(1)-x|+|\gamma(t)-x|\}\geq\frac{1+\sqrt{1+4\rho^2}}{2}\cdot\Delta,\]
and the proof is complete.
\end{enumerate}
\end{proof}

For the rest of the proof of Theorem \ref{thmAB'}, we set $\rho=\sqrt{8\alpha}$, and we denote by $\beta(\alpha)$ the inverse geometric mean of the two lower bounds in Proposition \ref{propkappa}:
\begin{equation}\label{betasubalpha}
\frac{1}{\beta(\alpha)}=\sqrt{\frac{1-\alpha}{(1+\alpha)\alpha}\cdot\frac{\left.\left(1+\sqrt{1+32\alpha}\right)\middle/2\right.-(4+\alpha)\alpha}{(1+\alpha)\alpha}}.
\end{equation}
With this definition, we have 
\[\frac{1}{\beta(\alpha)}=\frac{\sqrt{1+\alpha+o(\alpha)}}{\alpha}\quad\text{as $\alpha\to0^+$,}\]
which ensures that $\beta(\alpha)<\alpha$ for all sufficiently small $\alpha$.
This will be useful just below, see in particular the discussion following \eqref{eq:limsupfini}.

\subsubsection{Core of the proof}

\begin{proof}[Proof of Theorem \ref{thmAB'}]
Let $\mathcal{F}$ be a random closed subset of $[0,1]^d$, and assume that it satisfies Property $(\varnothing)$ with constants $\zeta>1$, and $Q>0$ and $q\in{]0,1[}$.
Recall that we want to prove the existence of a constant $s>1$ such that almost surely, for every non-trivial connected closed subset $\mathcal{C}\subset\mathcal{F}$, we have $\dim_H\mathcal{C}\geq s>1$.
To this end, fix a realisation of $\mathcal{F}$, and let $\mathcal{C}$ be a non-trivial connected closed subset of $\mathcal{F}$.
By the so-called energy method (see, e.g, \cite[Theorem 6.4.6]{bishop_peres_2016}), it suffices to construct  a Borel probability measure $\nu$ supported on $\mathcal{C}$ such that 
\begin{equation}\label{eqqenergy}
\int\int\frac{\mathrm{d}\nu(y)}{|x-y|^s}\mathrm{d}\nu(x)<\infty.
\end{equation}
In this direction, we first claim that it is possible to find a sequence $\left(\gamma_n:[0,1]\rightarrow\mathbb{R}^d\right)_{n\in\mathbb{N}^*}$ of paths, with ${\Delta:=\inf_{n\geq1}|\gamma_n[0,1]|>0}$, where $|A|$ denotes the diameter of $A\subset\mathbb{R}^d$, such that:
\begin{equation}\label{eqapproxcurve}
\text{for each $n\in\mathbb{N}^*$, we have $\gamma_n[0,1]\subset(\mathcal{C})_{1/n}$,}
\end{equation}
where $(\mathcal{C})_\varepsilon=\left\{x\in\mathbb{R}^d:d(x,\mathcal{C})\leq\varepsilon\right\}$ denotes the $\varepsilon$-neighborhood of $\mathcal{C}$.
Indeed, fix ${a\neq b\in\mathcal{C}}$.
For each $n\in\mathbb{N}^*$, the points $a$ and $b$ belong to the same connected component of 
\[\mathcal{O}_n=\left\{x\in\mathbb{R}^d:d(x,\mathcal{C})<1/n\right\}.\]
Since $\mathcal{O}_n$ is open, any connected component of $\mathcal{O}_n$ is path-connected, hence there exists a continuous path $\gamma_n:[0,1]\rightarrow\mathcal{O}_n$ that connects $a$ to $b$.
In particular, we have $\gamma_n[0,1]\subset(\mathcal{C})_{1/n}$, and the diameter $|\gamma_n[0,1]|$ of $\gamma_n[0,1]$ is at least $|a-b|>0$.
Following Aizenman and Burchard, we now use the splitting procedure $(\mathrm{P}_\alpha)$ recursively on the path $\gamma_n$, and derive a collection $\left(\mu^n_l\right)_{l\in\mathbb{N}}$ of Borel probability measures supported on $\gamma_n[0,1]$. 
Then, a deterministic property, which will be seen to hold for almost every realization of $\mathcal{F}$ thanks to Property $(\varnothing)$, will guarantee that it is possible to extract a sequence $\left(\nu_n=\mu^n_{L_n}\right)_{n\in\mathbb{N}^*}$, of which any subsequential weak limit $\nu$ will be a Borel probability measure supported on $\mathcal{C}$ such that \eqref{eqqenergy} holds.

More precisely now.
Fix a parameter $\alpha\in{\left]0,d^{-1/2}\right[}$ to be adjusted throughout the proof; for the moment, assume that $\alpha$ is small enough so that $\beta(\alpha)<\alpha$, where $\beta(\alpha)$ is defined in \eqref{betasubalpha}.
Given $\alpha$, fix $s=s(\alpha)>1$ such that $\beta(\alpha)<\alpha^s$.
Set $\delta_k=\alpha^k$ for all $k\in\mathbb{N}$, and denote by $k_0=k_0(\omega)$ the smallest integer $k\in\mathbb{N}$ such that $\delta_k\leq\Delta$. 
For each $n\in\mathbb{N}^*$, we split the path $\gamma_n$ into a collection $(\gamma^n_u,\,u\in\mathbb{T}^n)$ of subpaths, indexed by a plane tree $\mathbb{T}^n$ with root denoted by $o$, as follows.
First, by the definition of $\Delta$ and $k_0$, we have $|\gamma_n[0,1]|\geq\delta_{k_0}$.
Thus, there exists $t<t'\in[0,1]$ such that $|\gamma_n(t)-\gamma_n(t')|=\delta_{k_0}$, and we let $\gamma^n_o$ be the path $\theta\in[0,1]\mapsto\gamma_n((1-\theta)t+\theta t')$.
Then, by induction, having constructed the paths indexed by $\partial[\mathbb{T}^n]_l=\{u\in\mathbb{T}^n:|u|=l\}$, we apply for each $u\in\partial[\mathbb{T}^n]_l$ the procedure $(\mathrm{P}_\alpha)$ to the path $\gamma^n_u$, and we denote by $\gamma^n_{u1},\ldots,\gamma^n_{u\kappa_n(u)}$ the subpaths thereby generated. 
The children of $u$ in $\mathbb{T}^n$ are the nodes $u1,\ldots,u\kappa_n(u)$.
By construction, the following holds:
\begin{itemize}
\item for every $u\in\mathbb{T}^n$, we have ${|\gamma^n_u(0)-\gamma^n_u(1)|=\delta_{k_0+|u|}}$,
\item for any nodes $u,v$ that are not descendants of one another in $\mathbb{T}^n$, we have 
\[d(\gamma^n_u[0,1];\gamma^n_v[0,1])\geq\alpha\delta_{k_0+|u\wedge v|+1},\]
where $u\wedge v$ denotes the lowest common ancestor of $u$ and $v$.
\end{itemize}
Now, set $\pi_n(u)=\prod_{v\prec u}\kappa_n(u)^{-1}$ for all $u\in\mathbb{T}^n$, and let $\mu^n_l=\sum_{u\in\partial[\mathbb{T}^n]_l}\pi_n(u)\cdot(\gamma^n_u)_*\lambda$ for all $l\in\mathbb{N}$, where $(\gamma^n_u)_*\lambda$ denotes the push forward by $\gamma^n_u$ of the Lebesgue measure on $[0,1]$.
By construction, the measure $\mu^n_l$ is a probability supported on $\gamma_n[0,1]$, since $\sum_{u\in\partial[\mathbb{T}^n]_l}\pi_n(u)=1$ (this is easily checked by induction).

At this point, let us make the following calculation.
Let $\varepsilon_l:=\alpha\delta_{k_0+l+1}$ for all $l\in\mathbb{N}$.
Note that, since $\alpha<d^{-1/2}$ and $\Delta\leq\sqrt{d}$, we always have $\delta_{k_0}\geq\alpha\Delta$, and thus $\varepsilon_l\geq\alpha^3\Delta\cdot\alpha^l$ for all $l\in\mathbb{N}$.
For every $L\in\mathbb{N}$, we have
\begin{equation}\label{eqcalc}
\begin{split}
\int\int\frac{\mathrm{d}\mu^n_L(y)}{(\varepsilon_L\vee|x-y|)^s}\mathrm{d}\mu^n_L(x)&=\sum_{u,v\in\partial[\mathbb{T}^n]_L}\int_{\gamma^n_u[0,1]}\int_{\gamma^n_v[0,1]}\frac{\mathrm{d}\mu^n_L(y)}{(\varepsilon_L\vee|x-y|)^s}\mathrm{d}\mu^n_L(x)\\
&\leq\sum_{u,v\in\partial[\mathbb{T}^n]_L}\varepsilon_{|u\wedge v|}^{-s}\cdot\pi_n(u)\cdot\pi_n(v)\\
&=\sum_{l=0}^L\varepsilon_l^{-s}\cdot\sum_{\substack{u,v\in\partial[\mathbb{T}^n]_L\\|u\wedge v|=l}}\pi_n(u)\cdot\pi_n(v)\\
&\leq\sum_{l=0}^L\varepsilon_l^{-s}\cdot\sum_{t\in\partial[\mathbb{T}^n]_l}\sum_{\substack{u,v\in\partial[\mathbb{T}^n]_L\\u,v\succeq t}}\pi_n(u)\cdot\pi_n(v)\\
&=\sum_{l=0}^L\varepsilon_l^{-s}\cdot\sum_{t\in\partial[\mathbb{T}^n]_l}\pi_n(t)^2\\
&\leq\sum_{l=0}^L\varepsilon_l^{-s}\cdot\max_{t\in\partial[\mathbb{T}^n]_l}\pi_n(t)\cdot\sum_{t\in\partial[\mathbb{T}^n]_l}\pi_n(t)\\
&\leq\left(\alpha^3\Delta\right)^{-s}\cdot\sum_{l=0}^L\max_{t\in\partial[\mathbb{T}^n]_l}\pi_n(t)\cdot\alpha^{-sl}.
\end{split}
\end{equation}

Now, we claim that for almost every realization of $\mathcal{F}$, it is possible to choose $L=L_n(\omega)$, with $L_n\rightarrow\infty$ as $n\to\infty$, so that
\begin{equation}\label{eq:limsupfini}
\varlimsup_{n\to\infty}\sum_{l=0}^{L_n}\max_{t\in\partial[\mathbb{T}^n]_l}\pi_n(t)\cdot\alpha^{-sl}<\infty.
\end{equation}
It is here that the probabilistic machinery comes into play, and that we make use of the fact that $\mathcal{F}$ satisfies Property $(\varnothing)$. 
Suppose we introduce a family $(E_m)_{m\in\mathbb{N}^*}$ of events, with $\sum_{m\geq1}\mathbb{P}(E_m)<\infty$, such that $E_{k_0+l}$ is realized whenever there exists a node $u\in\partial[\mathbb{T}^n]_l$ with $\pi_n(u)>\beta(\alpha)^l$.
By the Borel-Cantelli lemma, almost surely, the event $E_m$ will fail to be realized for all sufficiently large $m$, which will prove \eqref{eq:limsupfini}: since $\beta(\alpha)<\alpha^s$, we have $\sum_{l\geq0}\beta(\alpha)^l\cdot\alpha^{-sl}<\infty$.

To get there, let $l\in\mathbb{N}^*$ and suppose that there exists $u\in\mathbb{T}^n_l$ such that $\pi_n(u)>\beta(\alpha)^l$.
Denoting by $u_0,\ldots,u_l$ the nodes on the geodesic path from the root to $u$ in $\mathbb{T}^n$, this can be reformulated as 
\[\prod_{0\leq k<l}\kappa_n(u_k)<\left(\frac{1}{\beta(\alpha)}\right)^l.\]
But by the definition of $\beta(\alpha)$ as the inverse geometric mean of the two lower bounds for $\kappa_n(\cdot)$ obtained in Proposition \ref{propkappa}, there must exist a number $j>l/2$ of indices ${l_1<\ldots<l_j\in\llbracket0,l\llbracket}$ such that, for each $i\in\llbracket1,j\rrbracket$, the path $\gamma^n_{u_{l_i}}$ does not deviate of a factor $\rho = \sqrt{8\alpha}$ from being a straight line.
In particular, there exists $\sigma_1\leq\ldots\leq\sigma_j<\tau_j\leq\ldots\leq\tau_1\in[0,1]$ such that, for every $i\in\llbracket1,j\rrbracket$:
\[|\gamma_n(\sigma_i)-\gamma_n(\tau_i)|=\delta_{k_0+l_i}\quad\text{and}\quad\gamma_n[\sigma_i,\tau_i]\subset S(\gamma_n(\sigma_i),\gamma_n(\tau_i);\rho\delta_{k_0+l_i}).\]
Now, writing $m_i=k_0+l_i$ for all $i\in\llbracket1,j\rrbracket$, let us discretise this information.

\medskip

\textsc{Discretisation step}\textbf{.}
For each $m\in\mathbb{N}$, let $\left(\overline{B}(z,\rho\delta_m)\right)_{z\in\mathcal{Z}_m}$ be a covering of $[0,1]^d$ by balls of radius $\rho\delta_m$, with centers $z\in[0,1]^d$ more than $\rho\delta_m$ apart so that the $\left(\overline{B}(z,\rho\delta_m/2)\right)_{z\in\mathcal{Z}_m}$ are disjoint.
For each ${i\in\llbracket1,j\rrbracket}$, we can find ${x_i,y_i\in\mathcal{Z}_{m_i}}$ such that ${\gamma_n(\sigma_i)\in\overline{B}(x_i,\rho\delta_{m_i})}$ and ${\gamma_n(\tau_i)\in\overline{B}(y_i,\rho\delta_{m_i})}$, and we have ${\gamma_n[\sigma_i,\tau_i]\subset S(x_i,y_i;2\rho\delta_{m_i})}$.
Discretising further, let us place a number $H\in\mathbb{N}^*$, to be adjusted soon, of points 
\[\left(z^i_h=\left(1-\frac{h}{H}\right)\cdot x_i+\frac{h}{H}\cdot y_i\right)_{h\in\llbracket0,H\rrbracket},\]
spread evenly on the line segment $[x_i,y_i]$.
By construction, the path $\gamma_n$ must meet each one of the balls $\left(\overline{B}\left(z^i_h,2\rho\delta_{m_i}\right)\,;\,h\in\llbracket0,H\rrbracket\right)$.
Now, since $\gamma_n[0,1]\subset(\mathcal{C})_{1/n}$, a similar statement holds for $\mathcal{F}$, namely: if $n\geq\left.1\middle/\left(\rho\delta_{m_j}\right)\right.$, then the set $\mathcal{F}$ must meet each one of the balls 
\[\left(\overline{B}\left(z^i_h,3\rho\delta_{m_i}\right)\,;\,\text{$i\in\llbracket1,j\rrbracket$ and $h\in\llbracket0,H\rrbracket$}\right).\]
Equivalently, the intersection event 
\[\text{$A^{m_1,\ldots,m_j}_{x_1,y_1;\ldots;x_j,y_j}$: ``for each $i\in\llbracket1,j\rrbracket$ and every $h\in\llbracket0,H\rrbracket$, the set $\mathcal{F}$ meets $\overline{B}\left(z^i_h,3\rho\delta_{m_i}\right)$''}\]
must be realized.
Here, the sequence $x_1,y_1;\ldots;x_j,y_j$ has the following properties:
\begin{itemize}
\item for every $i\in\llbracket1,j\rrbracket$, we have $x_i,y_i\in\mathcal{Z}_{m_i}$, with $(1-2\rho)\delta_{m_i}\leq|x_i-y_i|\leq(1+2\rho)\delta_{m_i}$,
\item for every $i\in\llbracket1,j\llbracket$, we have $x_{i+1},y_{i+1}\in S(x_i,y_i;2\rho\delta_{m_i})$.
\end{itemize}
We shall call any sequence satisfying those two properties \emph{admissible with respect to $m_1,\ldots,m_j$}. 

Summing up the previous reasoning, we have seen that, if there exists a node $u\in\partial[\mathbb{T}^n]_l$ such that ${\pi_n(u)>\beta(\alpha)^l}$, with $n\geq\left.1\middle/\left(\rho\delta_{k_0+l-1}\right)\right.$, then there must exist a number $j>l/2$ of indices ${l_1<\ldots<l_j\in\llbracket0,l\llbracket}$, and a sequence $x_1,y_1;\ldots;x_j,y_j$ which is admissible with respect to $k_0+l_1,\ldots,k_0+l_j$, such that the intersection event $A^{k_0+l_1,\ldots,k_0+l_j}_{x_1,y_1;\ldots;x_j,y_j}$ is realized. 
Let us now define, for all $m\in\mathbb{N}^*$, the event:
\[E_m=\bigcup_{m/3\leq j\leq m}\bigcup_{m_1<\ldots<m_j\in\llbracket0,m\llbracket}\bigcup_\text{$x_1,y_1;\ldots;x_j,y_j$ admissible}A^{m_1,\ldots,m_j}_{x_1,y_1;\ldots;x_j,y_j}.\]
If, for some $l\geq2k_0$ (this ensures that ${l/2\geq(k_0+l)/3}$), not too large so that $\rho\delta_{k_0+l-1}\geq1/n$, there exists a node $u\in\partial[\mathbb{T}^n]_l$ such that ${\pi_n(u)>\beta(\alpha)^l}$, then the event $E_{k_0+l}$ must be realized.
Now, let us show that upon adjusting the parameter $\alpha$, we can ensure that ${\sum_{m\geq1}\mathbb{P}(E_m)<\infty}$.

\medskip

\textsc{Summability of the $\mathbb{P}(E_m)$}\textbf{.}
Let $m\in\mathbb{N}^*$.
By the union bound, we have
\begin{equation}\label{equnionbound}
\mathbb{P}(E_m)\leq\sum_{m/3\leq j\leq m}\sum_{m_1<\ldots<m_j\in\llbracket0,m\llbracket}\sum_\text{$x_1,y_1;\ldots;x_j,y_j$ admissible}\mathbb{P}\left(A^{m_1,\ldots,m_j}_{x_1,y_1;\ldots;x_j,y_j}\right).
\end{equation}
Now, fix an integer $j$ such that $m/3\leq j\leq m$, fix indices $m_1<\ldots<m_j\in\llbracket0,m\llbracket$, and let $x_1,y_1;\ldots,x_j,y_j$ be an admissible sequence with respect to $m_1,\ldots,m_j$.
We will control the probability of the intersection event $A^{m_1,\ldots,m_j}_{x_1,y_1;\ldots;x_j,y_j}$ by making use of the fact that $\mathcal{F}$ satisfies Property $(\varnothing)$.
To this end, let us extract a collection of $\zeta$-separated balls from the 
\[\left(\overline{B}\left(z^i_h,3\rho\delta_{m_i}\right)\,;\,\text{$i\in\llbracket1,j\rrbracket$ and $h\in\llbracket0,H\rrbracket$}\right).\]
At this point, we choose ${H=\lfloor(1-2\rho)/((6\zeta+1)\rho)\rfloor}$. 
We have $|x_i-y_i|/H\geq(6\zeta+1)\rho\delta_{m_i}$ for each $i\in\llbracket1,j\rrbracket$, hence for any $h_1\neq h_2\in\llbracket0,H\rrbracket$, the following holds:
\begin{equation}\label{eqsep}
\text{for every $z_1\in\overline{B}\left(z^i_{h_1},3\zeta\rho\delta_{m_i}\right)$ and $z_2\in\overline{B}\left(z^i_{h_2},3\zeta\rho\delta_{m_i}\right)$, we have $|z_1-z_2|\geq\rho\delta_{m_i}$}.
\end{equation}
In particular, the balls 
\[\left(\overline{B}\left(z^j_h,3\rho\delta_{m_j}\right)\,;\,h\in\llbracket0,H\rrbracket\right)\]
are $\zeta$-separated: let us add them all to our collection.
To continue, note that the dilated balls 
\[\left(B\left(z^j_h,\zeta\cdot3\rho\delta_{m_j}\right)\,;\,h\in\llbracket0,H\rrbracket\right)\]
are all included in the sausage $S\left(x_j,y_j;3\zeta\rho\delta_{m_j}\right)$, which has diameter at most $(1+(6\zeta+2)\rho)\delta_{m_j}$ \big(indeed, the diameter of $S(x,y;r)$ is $|x-y|+2r$, and $|x_j-y_j|\leq(1+2\rho)\delta_{m_j}$\big).
Therefore, assuming now that $\alpha$ is small enough so as to have $(1+(6\zeta+2)\rho)<\rho/\alpha$: by \eqref{eqsep}, the sausage $S\left(x_j,y_j;3\zeta\rho\delta_{m_j}\right)$ meets the dilated ball 
\[\overline{B}\left(z^{j-1}_h,\zeta\cdot3\rho\delta_{m_{j-1}}\right)\]
for at most one $h_0\in\llbracket0,H\rrbracket$.
We add all the balls
\[\left(\overline{B}\left(z^{j-1}_h,3\rho\delta_{m_{j-1}}\right)\,;\,h\in\llbracket0,H\rrbracket\setminus\{h_0\}\right)\]
to our collection.
We iterate this argument, noticing that, as the sausages 
\[(S(x_i,y_i;3\zeta\rho\delta_{m_i})\,;\,i\in\llbracket1,j\rrbracket)\]
are nested (indeed, without loss of generality we may assume that $\alpha$ is small enough so that ${2+3\zeta\alpha\leq3\zeta}$), we only have to worry about intersections with the previous sausage at each step.
At the end of the construction, we obtain a collection of $\zeta$-separated balls that $\mathcal{F}$ must meet on the event $A^{m_1,\ldots,m_j}_{x_1,y_1;\ldots;x_j,y_j}$, which has cardinality at least ${(H+1)+(j-1)\cdot H\geq Hj}$. Since $\mathcal{F}$ satisfies Property $(\varnothing)$, we deduce that
\[\mathbb{P}\left(A^{m_1,\ldots,m_j}_{x_1,y_1;\ldots;x_j,y_j}\right)\leq Q\cdot q^{Hj}.\]
Coming back to \eqref{equnionbound}, we get
\[\begin{split}
\mathbb{P}(E_m)&\leq\sum_{m/3\leq j\leq m}\sum_{m_1<\ldots<m_j\in\llbracket0,m\llbracket}\sum_\text{$x_1,y_1;\ldots;x_j,y_j$ admissible}Q\cdot q^{Hj}\\
&=\sum_{m/3\leq j\leq m}\sum_{m_1<\ldots<m_j\in\llbracket0,m\llbracket}\#\{\text{admissible sequences with respect to $m_1,\ldots,m_j$}\}\cdot Q\cdot q^{Hj}.
\end{split}\]
Now, given an integer $j$ such that $m/3\leq j\leq m$, and indices $m_1<\ldots<m_j\in\llbracket0,m\llbracket$, let us control the number of admissible sequences with respect to $m_1,\ldots,m_j$.
First, there exists a constant $C=C(d,\alpha)>0$ such that $\#\mathcal{Z}_{m_1}\leq C\cdot\delta_{m_1}^{-d}$.
Indeed, the balls $\left(\overline{B}\left(z,\left.\rho\delta_{m_1}\middle/2\right.\right)\,;\,z\in\mathcal{Z}_{m_1}\right)$ are disjoint and included in the $\left.\rho\delta_{m_1}\middle/2\right.$-neighborood of $[0,1]^d$; thus, a volume argument yields:
\[\#\mathcal{Z}_{m_1}\leq\left(\frac{2\sqrt{d}}{\rho\delta_{m_1}}+1\right)^d\leq\left(\frac{2\sqrt{d}}{\rho}+1\right)^d\cdot\delta_{m_1}^{-d}=:C\cdot\delta_{m_1}^{-d}.\]
The constant $C$ depends on $d$ and $\rho=\sqrt{8\alpha}$, hence $C=C(d,\alpha)$.
Next, we claim that there exists a constant $c=c(d)$ such that for each $i\in\llbracket1,j\llbracket$,
\[\#\mathcal{Z}_{m_{i+1}}\cap S(x,y;2\rho\delta_{m_i})\leq\frac{c}{\rho}\cdot\left(\frac{\delta_{m_i}}{\delta_{m_{i+1}}}\right)^d\quad\text{for all $x,y\in\mathcal{Z}_{m_i}$ such that $|x-y|\leq(1+2\rho)\delta_{m_i}$.}\]
This is again by a volume argument, since the balls $\left(\overline{B}\left(z,\left.\rho\delta_{m_{i+1}}\middle/2\right.\right)\,;\,z\in\mathcal{Z}_{m_{i+1}}\cap S(x,y;2\rho\delta_{m_i})\right)$ are disjoint and included in the sausage $S\left(x,y;2\rho\delta_{m_i}+\left.\rho\delta_{m_{i+1}}\middle/2\right.\right)$.
Finally, we obtain that the number of admissible sequences with respect to $m_1,\ldots,m_j$ is bounded from above by
\[\left(C\cdot\delta_{m_1}^{-d}\right)^2\cdot\prod_{i=1}^{j-1}\left(\frac{c}{\rho}\cdot\left(\frac{\delta_{m_i}}{\delta_{m_{i+1}}}\right)^d\right)^2=C^2\cdot\frac{(c/\rho)^{2(j-1)}}{\delta_{m_j}^{2d}}=\frac{C^2}{(c/\rho)^2}\cdot\frac{(c/\rho)^{2j}}{\alpha^{2dm_j}}=:C'\cdot\frac{(c/\rho)^{2j}}{\alpha^{2dm_j}},\]
where $C'=C'(d,\alpha)$. 
Plugging this inequality into the above bound, we find
\[\begin{split}
\mathbb{P}(E_m)&\leq\sum_{m/3\leq j\leq m}\sum_{m_1<\ldots<m_j\in\llbracket0,m\llbracket}C'\cdot\frac{(c/\rho)^{2j}}{\alpha^{2dm_j}}\cdot Q\cdot q^{Hj}\\
&\leq\sum_{m/3\leq j\leq m}\sum_{m_1<\ldots<m_j\in\llbracket0,m\llbracket}C'\cdot\frac{(c/\rho)^{2m}\vee 1}{\alpha^{2dm}}\cdot Q\cdot q^{Hm/3}\\
&=C'\cdot Q\cdot\sum_{m/3\leq j\leq m}\binom{m}{j}\cdot\left(\frac{(c/\rho)^2 \vee 1}{\alpha^{2d}}\cdot q^{H/3}\right)^m\\
&\leq C'\cdot Q\cdot2^m\cdot\left(\frac{(c/\rho)^2 \vee 1}{\alpha^{2d}}\cdot q^{H/3}\right)^m=C'\cdot Q\cdot\left(2\cdot\frac{(c/\rho)^2 \vee 1}{\alpha^{2d}}\cdot q^{H/3}\right)^m.
\end{split}\]
Recalling that $\rho=\sqrt{8\alpha}$ and $H=\lfloor(1-2\rho)/((6\zeta+1)\rho)\rfloor$, a straightforward analysis shows that the term $2\cdot\left((c/\rho)^2 \vee 1\right) \cdot\alpha^{-2d}\cdot q^{H/3}$ can be made strictly smaller than $1$ by choosing $\alpha$ small enough.
For such $\alpha$, we get $\sum_{m\geq1}\mathbb{P}(E_m)<\infty$.

\medskip

\textsc{Concluding the proof}\textbf{.}
By the Borel--Cantelli lemma, almost surely, the event $E_m$ fails to be realized for all sufficiently large $m$.
Therefore, to almost every realization of $\mathcal{F}$ corresponds some $l_0\geq2k_0$ such that $E_m$ fails to be realized for all $m\geq k_0+l_0$.
Now, define $L_n=L_n(\omega)$ as the largest integer $l\in\mathbb{N}^*$ such that ${\rho\delta_{k_0+l-1}\geq1/n}$ (note that $L_n$ is well defined for all sufficiently large $n$, and that $L_n\to\infty$ as $n\to\infty$), and let $\nu_n=\mu^n_{L_n}$. Recalling \eqref{eqcalc}, we have
\[\int\int\frac{\mathrm{d}\nu_n(y)}{(\varepsilon_{L_n}\vee|x-y|)^s}\mathrm{d}\nu_n(x)\leq\left(\alpha^3\Delta\right)^{-s}\cdot\sum_{l=0}^{L_n}\max_{u\in\partial[\mathbb{T}^n]_l}\pi_n(u)\cdot\alpha^{-sl}.\]
By all the above work, if for some $l\in\llbracket2k_0,L_n\rrbracket$, there exists a node $u\in\partial[\mathbb{T}^n]_l$ such that ${\pi_n(u)>\beta(\alpha)^l}$, then the event $E_{k_0+l}$ must be realized.
Now we can write, recalling that $\beta(\alpha)<\alpha^s$:
\[\sum_{l=0}^{L_n}\max_{u\in\partial[\mathbb{T}^n]_l}\pi_n(u)\cdot\alpha^{-sl}\leq\sum_{l=0}^{l_0-1}\alpha^{-sl}+\sum_{l=l_0}^{L_n}\beta(\alpha)^l\cdot\alpha^{-sl}\leq\sum_{l=0}^{l_0-1}\alpha^{-sl}+\sum_{l\geq l_0}\left(\frac{\beta(\alpha)}{\alpha^s}\right)^l<\infty.\]
This proves that
\[\varlimsup_{n\to\infty}\int\int\frac{\mathrm{d}\nu_n(y)}{(\varepsilon_{L_n}\vee|x-y|)^s}\mathrm{d}\nu_n(x)<\infty.\]

Since we are working on the compact space $[0,1]^d$, the sequence of probability measures $(\nu_n)_{n\in\mathbb{N}^*}$ is automatically tight: let $\nu$ be any subsequential weak limit of $(\nu_n)_{n\in\mathbb{N}^*}$.
For each $\varepsilon>0$, by the Portmanteau theorem, we have ${\nu((\mathcal{C})_\varepsilon)\geq\varlimsup_{n\to\infty}\nu_n((\mathcal{C})_\varepsilon)=1}$ (the last equality holds because the support of $\nu_n$ is included in $(\mathcal{C})_\varepsilon$ for all sufficiently large $n$, thanks to \eqref{eqapproxcurve}). 
Since $\mathcal{C}$ is closed, we deduce that the probability measure $\nu$ is supported on $\mathcal{C}$. 
Furthermore, since $\nu$ is the weak limit of some subsequence $(\nu_{n_k})_{k\in\mathbb{N}^*}$, we have 
\[\begin{split}
\int\int\frac{\mathrm{d}\nu(y)}{(\varepsilon\vee|x-y|)^s}\mathrm{d}\nu(x)&=\lim_{k\to\infty}\int\int\frac{\mathrm{d}\nu_{n_k}(y)}{(\varepsilon\vee|x-y|)^s}\mathrm{d}\nu_{n_k}(x)\\
&\leq\varlimsup_{n\to\infty}\int\int\frac{\mathrm{d}\nu_n(y)}{(\varepsilon\vee|x-y|)^s}\mathrm{d}\nu_n(x)\\
&\leq\varlimsup_{n\to\infty}\int\int\frac{\mathrm{d}\nu_n(y)}{(\varepsilon_{L_n}\vee|x-y|)^s}\mathrm{d}\nu_n(x)<\infty.
\end{split}\]
This last upper bound does not depend on $\varepsilon$, and thus letting $\varepsilon\to0^+$, we conclude by the monotone convergence theorem that the integrability condition \eqref{eqqenergy} holds, completing the proof of Theorem \ref{thmAB'}.
\end{proof}

\subsection{Lower bound on the Hausdorff dimension of $\mathcal{F}_\infty$}\label{sec:lowerboundF}

The upper bound on the Hausdorff dimension of $\mathcal{F}_\infty$ stated in Theorem \ref{thm:main} was established in Proposition \ref{prop:upperbound}, now we come to the lower bound.
First, we check that $\mathcal{F}_\infty$ satisfies Property $(\varnothing)$, and that in dimension $d=2$, almost surely $\mathcal{F}_\infty$ has a non-trivial connected component.
With Theorem \ref{thmAB'}, this directly yields the result for $d=2$, which then bootstraps to any dimension with a slicing lemma. 

\begin{proposition}\label{propH2}
The frontier $\mathcal{F}_\infty$ satisfies Property $(\varnothing)$: there exists constants $Q>0$ and $q\in{]0,1[}$ such that, for every collection $\left(\overline{B}(x_i,r_i)\,;\,i\in\llbracket1,n\rrbracket\right)$ of $7$-separated balls, we have
\[\mathbb{P}\left(\text{for each $i\in\llbracket1,n\rrbracket$, the frontier $\mathcal{F}_\infty$ meets $\overline{B}(x_i,r_i)$}\right)\leq Q\cdot q^n.\]
\end{proposition}
\begin{proof}
This result is a consequence of Lemma \ref{lem:isolated}. 
Let $\left(\overline{B}(x_i,r_i)\,;\,i\in\llbracket1,n\rrbracket\right)$ be a collection of $7$-separated balls, and fix a realization of the intersection event: ``for each $i\in\llbracket1,n\rrbracket$, the frontier $\mathcal{F}_\infty$ meets $\overline{B}(x_i,r_i)$''. 
Notice that the initial points $R_0$ and $B_0$ belong to at most two distinct balls, with indices say $i_R$ and $i_B$. 
For every other $i\in\llbracket1,n\rrbracket$, both $R_0$ and $B_0$ lie outside $\overline{B}(x_i,7r_i)$, hence the event $G_i:=G^{x_i}_{7r_i/6,7r_i}$ of Lemma \ref{lem:isolated} fails to be realized. 
Indeed, if $G_i$ were realized, then the ball $\overline{B}(x_i,7r_i/6)$ would be monochromatic at the end of the coloring, and $\mathcal{F}_\infty$ would not meet $\overline{B}(x_i,r_i)$.  
Therefore, we have
\begin{multline*}
\left(\text{for each $i\in\llbracket1,n\rrbracket$, the frontier $\mathcal{F}_\infty$ meets $\overline{B}(x_i,r_i)$}\right)\\
\subset\left(\text{for each $i\in\llbracket1,n\rrbracket\setminus\{i_R,i_B\}$, the event $G_i$ fails to be realized}\right).
\end{multline*}
Since the events $(G_i)_{i\in\llbracket1,n\rrbracket}$ are independent and have probability at least $p>0$, we conclude that
\[\mathbb{P}\left(\text{for each $i\in\llbracket1,n\rrbracket$, the frontier $\mathcal{F}_\infty$ meets $\overline{B}(x_i,r_i)$}\right)\leq(1-p)^{n-2}=(1-p)^{-2}\cdot(1-p)^n.\]
This is the desired upper bound.
\end{proof}
 
Before proving that $\mathcal{F}_\infty$ has a non-trivial connected component, we first consider the following proposition.

\begin{proposition}\label{prop:Fnontrivial} 
Almost surely, the discrete frontier $\mathcal{F}_n=\left\{x\in[0,1]^d:d(x,\mathcal{R}_n)=d(x,\mathcal{B}_n)\right\}$ converges to $\mathcal{F}_{\infty}$ as $n\to\infty$, for the Hausdorff distance.
\end{proposition}
\begin{proof}
By the definition of $\mathcal{R}_\infty$ and $\mathcal{B}_\infty$, we have $\mathcal{R}_n\rightarrow\mathcal{R}_\infty$ and $\mathcal{B}_n\rightarrow\mathcal{B}_\infty$ as $n\to\infty$, for the Hausdorff distance (see \cite[Exercise 7.3.5]{burago^2ivanov}).
Now, fix $\varepsilon>0$, and let us prove that the inclusions $\mathcal{F}_n\subset(\mathcal{F}_\infty)_\varepsilon$ and $\mathcal{F}_\infty\subset(\mathcal{F}_n)_\varepsilon$ hold for all sufficiently large $n$.
\begin{itemize}
\item Since almost surely the set $\{R_0,B_0,X_1,X_2,\ldots\}$ is dense in $[0,1]^d$, for all sufficiently large $n$ the following holds: for each $x\in[0,1]^d$, the ball $\overline{B}(x,\varepsilon)$ contains an element of 
\[\{R_0,B_0,X_1,\ldots,X_n\}=\mathcal{R}_n\cup\mathcal{B}_n.\]
Then, for each $x\in\mathcal{F}_n$, as $d(x,\mathcal{R}_n)=d(x,\mathcal{B}_n)$, the ball $\overline{B}(x,\varepsilon)$ must contain an element of $\mathcal{R}_n$ and an element of $\mathcal{B}_n$.
Since $\mathcal{R}_\infty\cap\overline{B}(x,\varepsilon)$ and $\mathcal{B}_\infty\cap\overline{B}(x,\varepsilon)$ are two non-empty closed subsets whose union forms the connected set $\overline{B}(x,\varepsilon)$, they cannot be disjoint, hence $d(x,\mathcal{F}_\infty)\leq\varepsilon$.
This proves that $\mathcal{F}_n\subset(\mathcal{F}_\infty)_\varepsilon$.
\item Conversely, by the convergence of $\mathcal{R}_n$ and $\mathcal{B}_n$ towards $\mathcal{R}_\infty$ and $\mathcal{B}_\infty$, for all sufficiently large $n$ we have 
\[\sup_{x\in\mathcal{R}_\infty}d(x,\mathcal{R}_n)\leq\varepsilon\quad\text{and}\quad\sup_{x\in\mathcal{B}_\infty}d(x,\mathcal{B}_n)\leq\varepsilon.\]
Then, for every ${x\in\mathcal{F}_\infty=\mathcal{R}_\infty\cap\mathcal{B}_\infty}$, the ball $\overline{B}(x,\varepsilon)$ contains an element of $\mathcal{R}_n$ and an element of $\mathcal{B}_n$. 
Since 
\[\left\{y\in\overline{B}(x,\varepsilon):d(y,\mathcal{R}_n)\leq d(y,\mathcal{B}_n)\right\}\quad\text{and}\quad\left\{y\in\overline{B}(x,\varepsilon):d(y,\mathcal{R}_n)\geq d(y,\mathcal{B}_n)\right\}\]
are two non-empty closed subsets whose union forms the connected set $\overline{B}(x,\varepsilon)$, they cannot be disjoint; hence $d(x,\mathcal{F}_n)\leq\varepsilon$.
This proves that $\mathcal{F}_\infty\subset(\mathcal{F}_n)_\varepsilon$.
\end{itemize}
\end{proof}

\begin{corollary}\label{cor:comp} 
In dimension $d=2$, almost surely, the frontier $\mathcal{F}_{\infty}$ has a non-trivial connected component. 
\end{corollary}
\begin{proof}
Note that in dimension 2, the following holds: for each $n\in\mathbb{N}$, the discrete frontier $\mathcal{F}_{n}$ is a finite union of curves, where each curve is composed of line segments belonging to the boundary of the Voronoi cells of $R_0,B_0,X_1,\ldots,X_n$ (see Figure \ref{figvoronoi} below). 
By virtue of Lemma \ref{lem:isolated}, almost surely, there exists some $r>0$ such that the balls  $\overline{B}(R_0,r)$ and $\overline{B}(B_0,r)$ are monochromatic at the end of the coloring (see Remark \ref{remkeylemma}). 
In particular, this implies that the union of the red (resp.~blue) cells in the Voronoi diagram of the points $R_0,B_0,X_1,\ldots,X_n$ contains the ball $\overline{B}(R_0,r/2)$ \big(resp.~$\overline{B}(B_0,r/2)$\big). 
Therefore, the discrete frontier $\mathcal{F}_n$ contains a curve, i.e, the image of a continuous path $\gamma_n:[0,1]\rightarrow[0,1]^2$, of diameter at least $r$ (Figure \ref{figvoronoi} does not lie). 

\begin{figure}[ht]
\centering
\includegraphics[width=0.6\linewidth]{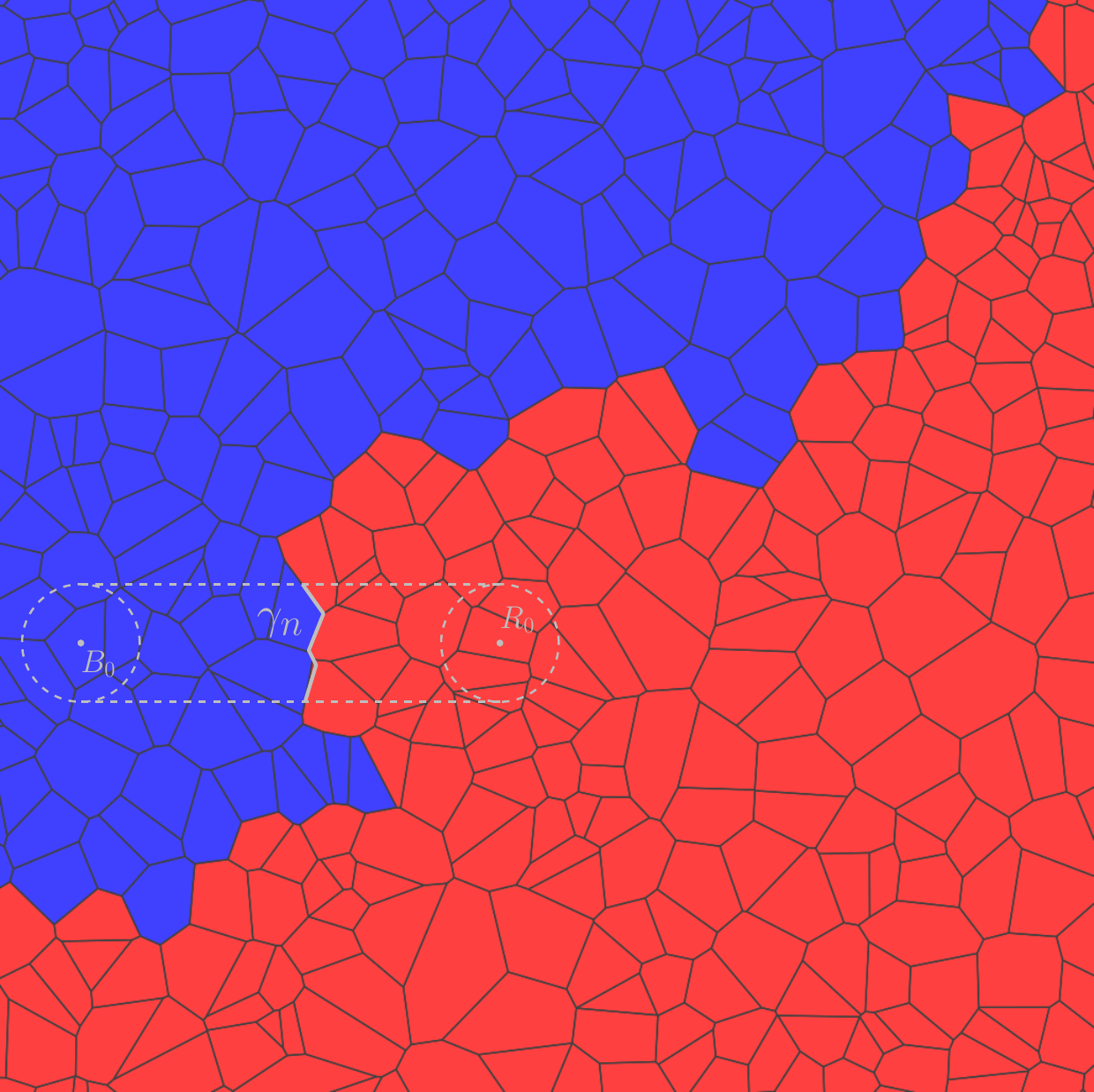}
\caption{Colored Voronoi diagram of the points $R_0,B_0,X_1,\ldots,X_n$. Balls with radius $r/2$ are represented around $R_0$ and $B_0$.
The path $\gamma_n$ is also represented, and it has diameter at least $r$.}\label{figvoronoi}
\end{figure}

Finally, recall that for the Hausdorff distance:
\begin{enumerate}[label=(\roman*)]
\item the set of closed subsets of $[0,1]^2$ is compact (see, e.g, \cite[Theorem 7.3.8]{burago^2ivanov}),
\item any limit of connected closed subsets is connected (see, e.g, \cite[Theorem 3.18]{falconer}).
\end{enumerate}
Therefore, we  can extract from $(\gamma_n[0,1])_{n\in\mathbb{N}}$ a subsequence which converges to some closed subset $\mathcal{C}\subset[0,1]^2$, that is necessarily connected, and has diameter at least $r>0$ (the diameter function is continuous for the Hausdorff distance).
Finally,  using that $\mathcal{F}_n$ converges to $\mathcal{F}_\infty$, the fact that  $\gamma_n[0,1]\subset \mathcal{F}_n$ implies that $\mathcal{C}$ is included in  $\mathcal{F}_\infty$.
\end{proof}

We now have all the ingredients to prove the lower bound of Theorem \ref{thm:main}.

\begin{proposition}\label{prop:lowerbound}
For any dimension $d \geq 2$, there exists $\varepsilon>0$ such that almost surely,
\[\dim_H\mathcal{F}_\infty\geq d-1+\varepsilon>d-1.\]
\end{proposition}
\begin{proof}
Fix $d\geq2$.
By Proposition \ref{propH2}, the frontier $\mathcal{F}_\infty$ satisfies Property $(\varnothing)$.
Therefore, by Theorem \ref{thmAB'}, there exists a constant $s>1$ such that almost surely, for every non-trivial connected closed subset $\mathcal{C}\subset\mathcal{F}_\infty$, we have
\[\dim_H\mathcal{C}\geq s>1.\]
If $d=2$, then by virtue of Corollary \ref{cor:comp}, almost surely such a subset $\mathcal{C}\subset\mathcal{F}_\infty$ exists, and we directly conclude that
\[\dim_H\mathcal{F}_\infty\geq s>1\quad\text{almost surely.}\]
For general $d\geq2$, we use a slicing argument.
Fix $\varepsilon\in{]0,s-1[}$.
Fix a plane $\mathcal{P}\subset\mathbb{R}^d$ that contains the red and blue seeds $R_0$ and $B_0$, and denote by $\mathcal{P}^\perp$ its orthogonal. 
By \cite[Theorem 1.6.2]{bishop_peres_2016}, there exists a constant $c=c(d)>0$ such that
\[H^{d-1+\varepsilon}(\mathcal{F}_\infty)\geq c\int_{\mathcal{P}^{\perp}}H^{1+\varepsilon}(\mathcal{F}_\infty\cap(z+\mathcal{P}))\mathrm{d}z,\]
where $H^t(A)$ denotes the $t$-dimensional Hausdorff measure of $A\subset\mathbb{R}^d$.
By Lemma \ref{lem:isolated}, almost surely there exists some $r>0$ such that the balls  $\overline{B}(R_0,r)$ and $\overline{B}(B_0,r)$ are monochromatic at the end of the coloring (see Remark \ref{remkeylemma}).
As in the proof of Corollary \ref{cor:comp}, this implies that for all $n\in\mathbb{N}$, the union of the red (resp.~blue) cells in the Voronoi diagram of $R_0,B_0,X_1,\ldots,X_n$ contains the ball  $\overline{B}(R_0,r/2)$ \big(resp.~$\overline{B}(B_0,r/2)$\big). Hence, for every $|z|\leq r/4$, the intersection between this union of red (resp.~blue) cells and the affine plane $(z+\mathcal{P})$
contains a two-dimensional  ball of radius $r/4$. Now, using the exact same arguments as in the proof of Corollary \ref{cor:comp}, we deduce that for every $|z|\leq r/4$, the closed subset $\mathcal{F}_\infty\cap(z+\mathcal{P})\subset\mathcal{F}_\infty$ has a non-trivial connected component. Therefore, almost surely, we have $\dim_H(\mathcal{F}_\infty\cap(z+\mathcal{P}))\geq s$.
By the definition of $\varepsilon$, it follows that almost surely, we have $H^{1+\varepsilon}\left(\mathcal{F}_\infty\cap(z+\mathcal{P})\right) = \infty$ for all $|z|\leq r/4$, hence
\[H^{d-1+\varepsilon}(\mathcal{F}_\infty)\geq c\int_{\mathcal{P}^{\perp}\cap\overline{B}(0,r/4)}H^{1+\varepsilon}\left(\mathcal{F}_\infty\cap(z+\mathcal{P})\right)\mathrm{d}z=\infty.\]
This proves that $\dim_H\mathcal{F}_\infty\geq d-1+\varepsilon$ almost surely.
\end{proof}

\section*{Acknowledgements}
We warmly thank David Aldous for discussions about \cite{aldous2018random} and for providing us with the reference \cite{preater2009species}.
We are grateful to the participants of the PizzaMa seminar, during which this work was initiated.
Finally, we are thankful to the anonymous referee for their careful reading of the paper, and their various comments.

\bibliographystyle{alea3}
\bibliography{bibli.bib}

\end{document}